\documentclass[10pt]{amsart}

\usepackage{amssymb}
\usepackage{amsmath}
\newcommand{\proclaim}[2]{\medbreak {\bf #1}{\sl #2} \medbreak}

\let\newpf\proof \let\proof\relax 
\newenvironment{pf}{\newpf[\proofname]}{\qed\endtrivlist}

\def\ssubset{\Subset}

\def\top{{\mathrm{top}}}

\def\be{\begin{equation}}
\def\ee{\end{equation}}

\def\op{\overline \partial}

\def\ba{{\begin{align}}}
\def\ea{{\end{align}}}

\def\0{{\mathbf 0}}

\def\cal{\mathcal}

\newcommand{\Fix}{\operatorname{Fix}}

\newtheorem{thm}{Theorem}[section]
\newtheorem{cor}[thm]{Corollary}
\newtheorem{lem}[thm]{Lemma}
\newtheorem{lemma}[thm]{Lemma}

\newtheorem{prop}[thm]{Proposition}
\newtheorem{schw}{Schwarz Lemma} 

\theoremstyle{remark}
\newtheorem{rem}{Remark}[section]

\newtheorem{example}{Example}[section]

\numberwithin{equation}{section}

\def \bn {\hfill \\ \smallskip\noindent}

\theoremstyle{definition}

\def\proof{\bn {\bf Proof.} }

\def\note#1
{\marginpar
{\tiny $\leftarrow$
\par
\hfuzz=20pt \hbadness=9000 \hyphenpenalty=-100 \exhyphenpenalty=-100
\pretolerance=-1 \tolerance=9999 \doublehyphendemerits=-100000
\finalhyphendemerits=-100000 \baselineskip=6pt
#1}\hfuzz=1pt}

\newcommand{\bignote}[1]{\begin{quote} \sf #1 \end{quote}}

\newcommand{\di}{\partial}
\newcommand{\dibar}{\bar\partial}
\newcommand{\ra}{\rightarrow}

\def\ssk{\smallskip}
\def\msk{\medskip}

\def\nin{\noindent}

\def\sm{\smallsetminus}

\newcommand{\diam}{\operatorname{diam}}
\newcommand{\dist}{\operatorname{dist}}

\newcommand{\cl}{\operatorname{cl}}
\newcommand{\inter}{\operatorname{int}}
\renewcommand{\mod}{\operatorname{mod}}
\newcommand{\tl}{\tilde}

\newcommand{\orb}{\operatorname{orb}}

\newcommand{\id}{\operatorname{id}}

\newcommand{\Dil}{\operatorname{Dil}}

\newcommand{\isom}{\approx}

\renewcommand{\Im}{\operatorname{Im}\,}

\newcommand{\eps}{{\epsilon}}
\newcommand{\De}{{\Delta}}
\newcommand{\de}{{\delta}}
\newcommand{\la}{{\lambda}}

\newcommand{\si}{{\sigma}}
\newcommand{\Om}{{\Omega}}
\newcommand{\om}{{\omega}}

\renewcommand{\AA}{{\cal A}}
\newcommand{\BB}{{\cal B}}
\newcommand{\CC}{{\cal C}}

\newcommand{\EE}{{\cal E}}

\newcommand{\II}{{\cal I}}
\newcommand{\FF}{{\cal F}}
\newcommand{\GG}{{\cal G}}

\newcommand{\HH}{{\cal H}}
\newcommand{\KK}{{\cal K}}

\newcommand{\MM}{{\cal M}}

\newcommand{\QQ}{{\cal Q}}

\newcommand{\SSS}{{\cal S}}

\newcommand{\UU}{{\cal U}}
\newcommand{\VV}{{\cal V}}

\newcommand{\ZZ}{{\cal Z}}

\newcommand{\C}{{\mathbb C}}
\newcommand{\D}{{\mathbb D}}

\newcommand{\I}{{\mathbb I}}

\newcommand{\N}{{\mathbb N}}

\newcommand{\R}{{\mathbb R}}
\newcommand{\T}{{\mathbb T}}

\newcommand{\Z}{{\mathbb Z}}

\def\B0{{\bold{0}}}

\catcode`\@=12

\def\Empty{}
\newcommand\oplabel[1]{
  \def\OpArg{#1} \ifx \OpArg\Empty {} \else
  	\label{#1}
  \fi}
		
\newcommand{\comm}[1]{}
\newcommand{\comment}[1]{}

\begin{document}

\title[Exponential contraction of renormalization]
{The full renormalization horseshoe \\ 
for unimodal maps of higher degree: \\
\tiny{exponential contraction along hybrid classes}}

\date{\today}

\author{Artur Avila and Mikhail Lyubich}

\address{
CNRS UMR 7586,
Institut de Math\'ematiques de Jussieu\\
175 rue du Chevaleret\\
75013--Paris, France \&
IMPA\\
Estrada Dona Castorina, 110\\
22460-320, Rio de Janeiro, RJ, Brazil.
}

\email{artur@math.sunysb.edu}
\urladdr{www.impa.br/$\sim$avila}

\address{
Department of Mathematics, Stony Brook University\\
Stony Brook, NY\\
11794-3651 \&
Department of Mathematics, University of Toronto\\
Room 6290, 40 St. George Street, Toronto, ON M5S 2E4, Canada.
}

\email{mlyubich@math.sunysb.edu}

\def\IMSmarkvadjust{0 pt}
\def\IMSmarkhadjust{0 pt}
\def\IMSmarkhpadding{0 pt}
\def\IMSpubltext{Published in modified form:}
\def\SBIMSMark#1#2#3{
 \font\SBF=cmss10 at 10 true pt
 \font\SBI=cmssi10 at 10 true pt
 \setbox0=\hbox{\SBF \hbox to \IMSmarkhpadding{\relax}
                Stony Brook IMS Preprint \##1}
 \setbox2=\hbox to \wd0{\hfil \SBI #2}
 \setbox4=\hbox to \wd0{\hfil \SBI #3}
 \setbox6=\hbox to \wd0{\hss
             \vbox{\hsize=\wd0 \parskip=0pt \baselineskip=10 true pt
                   \copy0 \break%
                   \copy2 \break%
                   \copy4 \break}}
 \dimen0=\ht6   \advance\dimen0 by \vsize \advance\dimen0 by 8 true pt
                \advance\dimen0 by -\pagetotal
	        \advance\dimen0 by \IMSmarkvadjust
 \dimen2=\hsize \advance\dimen2 by .25 true in
	        \advance\dimen2 by \IMSmarkhadjust

%
%
  \openin2=publishd.tex
  \ifeof2\setbox0=\hbox to 0pt{}
  \else 
     \setbox0=\hbox to 3.1 true in{
                \vbox to \ht6{\hsize=3 true in \parskip=0pt  \noindent  
                {\SBI \IMSpubltext}\hfil\break
                \input publishd.tex 
                \vfill}}
  \fi
  \closein2
  \ht0=0pt \dp0=0pt
 \ht6=0pt \dp6=0pt
 \setbox8=\vbox to \dimen0{\vfill \hbox to \dimen2{\copy0 \hss \copy6}}
 \ht8=0pt \dp8=0pt \wd8=0pt
 \copy8
 \message{*** Stony Brook IMS Preprint #1, #2. #3 ***}
}

\SBIMSMark{2010/3}{May 2010}{}

\begin{abstract}

We prove exponential contraction of renormalization along hybrid classes of
infinitely renormalizable unimodal maps (with arbitrary
combinatorics), in any even degree $d$.  We then conclude that orbits of
renormalization are asymptotic to the full renormalization horseshoe, which
we construct.  Our argument for exponential contraction
is based on a precompactness property of the renormalization
operator (``beau bounds''), which is leveraged in the abstract
analysis of holomorphic iteration.
Besides greater generality, it yields a unified approach to all
combinatorics and degrees: there is no need to account for the varied
geometric details of the dynamics, which were the typical source of
contraction in previous restricted proofs.


\end{abstract}

\setcounter{tocdepth}{1}

\maketitle

\tableofcontents

\section{Introduction}

\subsection{Renormalization Conjecture and Regular or Stochastic Theorem}

The Renormalization Conjecture formulated in mid 1970's by Feigenbaum \cite{F} and independently by 
Coullet and Tresser \cite{TC} has been a focus of research ever since.
Roughly speaking, it says that a certain  ``renormalization  
is hyperbolic in an appropriate infinite-dimensional functional space. 
It explains remarkable universality properties on various families of
dynamical systems (see  \cite{C} for a collection of early papers on the
subject).  More recently, it has played a central role in the
measure-theoretical analysis of one-dimensional dynamical systems,
particularly in the proofs of the Regular or Stochastic
Dichotomy in the real quadratic family \cite{regular}
and more general spaces of quasiquadratic unimodal maps  \cite{ALM}.  

Here we will consider the renormalization operator $R$ in the
space $\CC_d^\R$  of real 
unicritical polynomial-like maps of an arbitrary even degree $d\geq 2$. 
Hyperbolicity of $R$ was proven
for bounded combinatorics
by Sullivan, McMullen and one of the authors in \cite{S,McM,universe},
and then for all combinatorics in the quadratic case \cite{regular}. 
Our goal is to  generalize the latter result to an arbitrary
even degree $d$.  

In this paper we will prove that the renormalization operator $R$ 
has an invariant horseshoe $\AA$  
and is exponentially contracting on the corresponding hybrid lamination. 
In the forthcoming paper we will deal with the transverse
unstable direction.  Together with the previous analysis of
non-renormalizable unimodal maps \cite{ALS}, 
this will prove the Regular or Stochastic Dichotomy
in any unicritical family $p_c: x\mapsto x^d+c$ ($d \geq 2$ even): 
For almost any real $c$ (for which $p_c$ has an invariant interval),
the map $p_c$ is either regular (i.e., it has an attracting cycle)
or stochastic (i.e., it has an absolutely continuous invariant measure).

Besides supplying a more general version of the Renormalization Theorem,
our goal is to address the issue of exponential contraction along hybrid
leaves in a novel unified way, which does not involve fine geometric
considerations (highly dependent on the combinatorics and degree).
Our approach simplifies the previously known proofs in the quadratic-like
case, even for the renormalization with bounded combinatorics. 
Namely, we will derive the desired result from the
previously known {\it beau bounds} for real maps \cite{S,LvS,LY},
and basic facts of functional analysis and topology.

\subsection{Statement of the result}

Let us now formulate our main result more precisely. 
To this end we need a few basic definitions that we now outline;
a more detailed background will be supplied in the main body of the paper. 

A {\it unicritical polynomial-like map} of degree $d$  is a degree
$d$ branched covering
$f: U\ra V$ between two topological disks $U\Subset V$ that has a single
critical point.  We normalize $f$
so that $f(z)= z^d+c+O(z^{d+1})$ at the origin.%
\footnote{Note that this normalization survives rotations through
$2\pi k/ (d-1)$, $k\in \Z/ (d-1)\Z$:  
conjugating a normalized map by such a rotation, 
we obtain a normalized map with rotated $c$.} 

The (filled) {\it Julia set} $K(f)$ is the set of non-escaping points.
It is either connected or a Cantor set depending on whether
$0\in K(f)$ or not.
If $f$ is a polynomial-like map with connected Julia set then the
corresponding polynomial-like {\it germ}
is defined as the class of polynomial-like maps $\tl f$ with the same
Julia set and such that $\tl f|\, K(f)= f|\, K(f)$.  

Let $\CC=\CC_d$ stand for the space of normalized polynomial-like germs
of degree $d$ with connected Julia set. It intersects the polynomial
family $p_c: z\mapsto z^d+c$, $c\in \C$,  by the {\it Multibrot set}
$\MM= \MM_d$ (defined as the set of $c$ for which the Julia set
$K(p_c)$ is connected).

Two polynomial-like germs are called {\it hybrid equivalent} if they have
representatives 
$f: U\ra V$ and $\tilde f: \tilde U\ra \tilde V$ that are conjugate 
by a quasiconformal homeomorphism $h: V\ra \tilde V$  such that
$\dibar h=0$ almost everywhere on $K(f)$. 
The corresponding equivalence classes are called {\it hybrid classes}. 
According to the Douady-Hubbard Straightening Theorem \cite{DH},   
any hybrid class in $\CC$ intersects the Multibrot set $\MM$ 
by  an orbit of the rotation group
$\Z / (d-1)\Z$.

A unicritical polynomial-like germ $f$ is called {\it renormalizable} 
if there is a disk $\Om\ni 0$ and a $p\geq 2$ such that the map
$f^p|\, \Om$ is unicritical  polynomial-like with connected Julia set
(subject of a a few extra technical requirements -- 
see \S \ref{renorm sec}). Appropriately normalizing this
polynomial-like germ,
we obtain the {\it renormalization} of $f$.  
If $p$ is the smallest period for which $f$ is renormalizable,
then the corresponding renormalization is denoted $Rf$. 

We can now naturally define {\it infinitely renormalizable
polynomial-like germs}.   
Let $\II^\R$ stand for the space of {\it real} infinitely renormalizable
polynomial-like germs (that is, the germs preserving the real line),
and let $\II^{(\R)}$ stand for the space of polynomial-like germs
that are hybrid equivalent to the real ones. 
The renormalization operator $R$ naturally acts in both spaces
preserving the hybrid partition.
In what follows, this partition will serve as
the {\it stable lamination}:  

\msk

\noindent{\bf Main Theorem}.

{\it  There is an $R$-invariant precompact set $\AA\subset \II^\R$
(the renormalization horseshoe) such that $R|\, \AA$ is topologically
conjugate to the two-sided shift in infinitely many symbols,
and any germ $f\in \II^{(\R)}$ is attracted to some orbit of
$\AA$ at a uniformly exponential rate, in a suitable ``Carath\'eodory metric''.}

\msk

See Theorem \ref {real horseshoe} for a slightly more detailed
formulation.

\begin{rem}

Our approach to exponential contraction also applies to certain non-real
renormalization combinatorics (for which the appropriate beau bounds have
been established, see \cite {K}, \cite {KL1}, \cite {KL2}).  See Theorem
\ref {contr}.

\end{rem}

\subsection{Outline of the proof}

We start with the argument for exponential contraction along hybrid
classes.  To fix ideas, 
let us consider first the case of a {\it fixed hybrid leaf} $\HH_c$ (the
connected component of the hybrid class of $p_c$)
so that every time we iterate the same renormalization operator $R$.

\msk {\it Hybrid lamination}.

In \S \ref{hybrid classe sec} we endow hybrid leaves with a
{\it path holomorphic structure} 
and show that all of them are bi-holomorphically equivalent. 
The path holomorphic structure endows these spaces with
{\it Carath\'eodory pseudo-metrics}
and we prove that they are {\it Carath\'eodory hyperbolic}, i.e., 
these pseudo-metrics are, in fact, metrics.

\msk {\it The Schwarz Lemma}. The renormalization
operator, as a map from one hybrid leaf to another,
is holomorphic with respect to their path holomorphic
structure.  This puts us in a position to apply the
Schwarz Lemma to the analysis of its iterates: 
its weak form (see \S \ref{Schwarz lemma sec})
already implies that renormalization is {\it weakly contracting}
with respect to the Carath\'eodory metric.

\msk

  {\it Beau bounds in $\HH_c$}
mean by definition that there exists a compact set $\KK
\subset \CC$ such that for every $f \in \HH_c$,
$R^n f \in \KK$ for any $n$ sufficiently large (depending only on the
quality of the analytic extension of $f$), see \S \ref{beau sec}.

For real maps with stationary combinatorics,
beau bounds were proved by Sullivan (see \cite{S,MvS}) in early 1990's,
for complex maps with {\it primitive} stationary combinatorics,
they have been recently established by Kahn \cite{K}.
(Note that this result covers all real stationary combinatorics except
the period doubling.)  Exponential contraction can be easily concluded from
beau bounds through the entire complex hybrid leaf $\HH_c$:

\comm{

\msk {\it Rigidity}.

Let us assume that we have beau bounds for the stationary combinatorics
in question. Then there exists only one (up to rotation) infinitely
renormalizable polynomial
$p_c: z\mapsto z^d+c$ with this combinatorics (see \cite{S} for real maps and
\cite{puzzle} for complex ones). 
Then  one of the spaces $\HH_c$ is {\it invariant} under the
renormalization $R$, and the latter is path holomorphic. By the Schwarz
Lemma (in a weak form), it is automatically {\it weakly contracting}
with respect to the Carath\'eodory metric.   
}

\msk {\it The Strong Schwarz Lemma} shows that holomorphic endomorphisms are
strongly contracting with respect to the Carath\'eodory metric, provided
the image is ``small'' in the range, in the sense that the diameter is less
than $1$,\footnote{We gauge the Carath\'eodory metric so that any space has diameter
at most 1.  Condition that the diameter of $\KK$ in $\QQ$ is less than
1 means that ``$\KK$ is well inside of $\QQ$''.}
see \S \ref{Schwarz lemma sec}.
Beau bounds imply that for any compact set $\QQ \subset \HH_c$ there exists
$N$ such that $R^n(\QQ)$ contained in the (universal) compact set $\KK$ for
$n \geq N$.
By selecting
$\QQ\supset \KK$ sufficiently large, we fulfill the smallness condition of
$\KK$ inside $\QQ$, and
can conclude that
     $R^N|\, \QQ$ is strongly contracting (\S \ref{beau sec}).

\msk 

We will now give a different argument for stationary real combinatorics
that relies only on the beau bounds for {\it real} maps.
It makes use of one general idea of functional analysis:

\msk

{\it Almost periodicity}.  The beau bounds for real maps (and even just for
real polynomials) imply that the cyclic semigroup $\{R^n\}_{n=0}^\infty$ 
is precompact in the topology of uniform convergence
on compact sets of $\HH_c$.
(Such semigroups are called almost periodic, see \cite{Lyu,Lyu2}.)
Then the $\omega$-limit set of this semigroup is a {\it group}.
Its unit element is  a retraction $P: \HH_c\ra \ZZ= \Fix P$.

\begin{rem}

More directly, the precompactness of $\{R^n\}$ allows us to find nearby iterates
$R^n$ and $R^m$ with $m>2n$.  It follows that $R^{n-m}$ is close to the
identity in $\mathrm {Im} R^n \supset \mathrm {Im} R^{n-m}$.  Taking limits we get
a map $P$ which is exactly the identity in $\mathrm {Im} P$, so that $P$ is a
retraction.
\end{rem}

Let $P^\R: \HH_c^\R\ra \ZZ^\R$ be the restriction of $P$ to the real slice.

\msk

{\it Topological argument and analytic continuation)} (\S \ref{triviality sec}).
The beau bounds for real maps imply that the real slice $\ZZ^\R$ is
compact.  By the Implicit Function Theorem, $\ZZ^\R$ is a
finite-dimensional manifold.
But one can show that the space $\HH_c^\R$ is contractible
(see Lemma \ref{contractible} and Theorem \ref {productstructure}), and
hence the retract $\ZZ^\R$ is contractible as well. But the only
contractible compact finite dimensional manifold
(without boundary) is a single point. 


So, $P^\R$ collapses the real slice  $\HH_c^\R$ to a single point $f_*$.
Since $P$ is holomorphic, it collapses the whole space $\HH_c$ to 
$f_*$ as well.

\msk


Since $P$ is constant it follows that $R^n\to P$ uniformly on compact subsets
of $\HH_c$, and we can conclude exponential contraction through the strong
Schwarz Lemma as before.

This completes the argument for the case of stationary combinatorics.

\bigskip {\it Unbounded combinatorics.}
Beau bounds for arbitrary real maps were established in \cite{LvS,LY}.
For complex maps, they have been recently established for a
fairly big class of combinatorics
in \cite{KL1,KL2} (see also \cite{puzzle} for earlier results).
Our first argument that uses the beau bounds for complex maps extends
to the unbounded combinatorics case in 
a straightforward way, using the leafwise Carath\'eodory metric on the
hybrid lamination (\S \ref{beau sec}). 
The second argument based on almost periodicity requires an extension of
this idea from semigroups to {\it cocycles} (or {\it grupoids}).
This is carried out in \S \ref{cocycles}.

\bigskip {\it Horseshoe.}

Once contraction is proved, the horseshoe is constructed in the familiar way (see \cite{regular})
using {\it rigidity} of real maps \cite{puzzle,GS,KSS}.  It is
automatically semi-conjugate to the full shift.  A further argument
based on the analysis of the analytic continuation of anti-renormalizable
maps (which becomes substantially more involved in the higher degree
case, see Appendix \ref{analytic}) yields the full topological conjugacy.

\subsection{Comparison with earlier approaches}

In the case of stationary combinatorics, two approaches were previously
used to construct the fixed point $f_*$
and to prove convergence to $f_*$ in the hybrid class $\HH(f_*)$. 
The first one, due to Sullivan, is based on ideas
of Teichm\"uller theory (see \cite{S,MvS});
the other one, due to McMullen, is based on
a geometric theory  of {\it towers} and their
quasiconformal rigidity  \cite{towers}.

The Teichm\"uller approach, albeit beautiful and natural,
faces a number of subtle technical issues. 
Also, it does not seem to lead to the exponential contraction.\footnote
{It has been suggested that this relates to the fact the Teichm\"uller
approach naturally deals with conformal (rather than affine) equivalence
between the polynomial-like germs.  However, our Schwarz Lemma argument
seems to work equally well for conformal classes.}
The geometric tower approach can be carried all the way to prove
exponential contraction \cite{towers}.
On the other hand, in \cite{universe}, exponential contraction
was obtained by combining towers rigidity (as a source of
contraction, but without the rate) with the Schwarz Lemma in 
Banach spaces. 

Both approaches generalize without problem to the bounded
combinatorics case.
The tower approach can be carried further to ``essentially bounded''
combinatorics \cite{Hi}. 
The  remaining ``high'' combinatorics case
(as well as the oscillating situation)
 was handled in \cite{regular}
using the geometric property of growth of moduli in the Principal Nest
of the  Yoccoz puzzle \cite{puzzle}.
This is a powerful geometric property which is valid only in the
quadratic case.
So, this method is not sufficient in  the higher degree case
(at least, it would require further non-trivial geometric analysis). 

The approaches developed in this paper use much softer
geometric input (only beau bounds)
and treat all the combinatorial cases in a unified way. 

\begin{rem}

  There are also computer-assisted methods going back
to the classical paper by Lanford \cite{La},
as well as approaches that do not rely on
holomorphic dynamics \cite{E,Ma}.
These methods can be important for dealing with the case of
fractionary degree $d$.
The almost periodicity idea can possibly contribute to it, too.

\end{rem}

\subsection{Basic Notation}


$\D = \{ z: |z|<1\}$ is the unit disk, \\
$\D_r= \{ z: |z|< r \}$ is the disk of radius $r$,\\ 
and $\T=\{ z: |z|=1\}$ is the  unit circle;\\
$p_c: z\mapsto z^d+c$ is the unicritical polynomial family.  

We assume the reader's familiarity with the basic theory of quasiconformal (``qc'') maps.
We let $\Dil h$ be the dilatation of a qc map.

{\bf Acknowledgements:}  The first author would like to thank the
hospitality of the University of Toronto, the Fields Institute and the
Stony Brook University.  This research was partially conducted during the
period the first author served as a Clay Research Fellow.
The second author was partially supported by the NSF and NSERC.

\def\Lip{\mathrm{Lip}}
\def\Hol{\mathrm{Hol}}
\def\MI{\mathcal{MI}}

\section{Hybrid classes, external maps, and renormalization}
Theory of polynomial-like maps was laid down in \cite{DH}
and further developed, particularly in the quadratic-like setting,  in  \cite{McM,universe}. 
In this section we will refine the basic theory  
in the case of  unicritical polynomial-like maps of arbitrary degree.

\subsection{Holomorphic motions} \label {hol motions}                   
Given a domain $D\subset \C$ with a base point $\la_0$
and a set  $X_0\subset\C$, a {\it holomorphic motion} of $X_0$ over $D$ is a family
of injections $h_\la: X_0\ra\C$, $\la\in D$, such that $h_{\la_0}=\id$ and $h_\la(z)$ is
holomorphic in $\la$ for any $z\in X_0$. Let $X_\la = h_\la(X_0)$.

We will summarize
fundamental properties of holomorphic motions which are usually referred to as the
$\la$-{\it lemma}. It consists of two  parts: extension of the motion and transverse
quasiconformality, which will be stated separately.

\proclaim {Extension $\la$-Lemma \cite{Sl}.}
{
 A holomorphic motion $h_\la: X_0\ra X_\la$  of a set
$X_0\subset\C$ over the disk $\D$ 
admits an extension to a holomorphic motion
$\hat h_\la: \C\ra\C$ of the whole complex plane over $\D$.  
}

\begin{rem}
   We will usually keep the same notation, $h_\la$,  for the extended motion.
\end{rem}

\proclaim {Quasiconformality $\la$-Lemma \cite{MSS}.}
{
Let $h_\la: U_0\ra U_\la$ be a holomorphic motion of a domain $U_0\subset\C$
over the disk $\D$, based on $0$. Then over any smaller disk $\D_r$, $r<1$, all the maps
$h_\la$ are $K(r)$-qc, where $K(r)=\frac {1+r} {1-r}$.}

\subsection{Polynomial-like maps and germs}
A  {\it polynomial-like map} (``p-l map'') of degree $d\geq 2$
is a holomorphic branched covering $f:U \to V$ of degree $d$
between quasidisks $U \Subset V$.  
Its filled Julia set is $K(f)=\cap_{n \geq 0} f^{-n}(U)$, 
and the Julia set is $J(f)=\partial K(f)$.
In what follows, the letter $d$ will be reserved for the degree of $f$.

A p-l map is called {\it unicritical} if it has a unique critical point
(of local degree $d$). We will normalize unicritical polynomial-like maps 
so that $0 \in U$ is the critical point and
$f(z)=z^d+c+O(z^{d+1})$ near $0$.
In what follows,  polynomial-like maps under consideration will be
assumed unicritical.  The annulus $V\sm U$ is called the
{\it fundamental annulus} of a p-l map $f: U\ra V$
(the corresponding open and closed annuli, $V\sm \bar U$ and
$\bar V\sm U$ will also
be called ``fundamental'').

Basic examples of p-l maps are provided by appropriate restrictions of
unicritical polynomials $p_c:z \mapsto z^d+c$, e.g., $p_c:\D_r \to
p_c(\D_r)$ for $r>1+|c|$.

The {\it Basic Dichotomy} asserts that the (filled) Julia
set of $f$ is either connected or a Cantor set,
and the former happens iff $0\in K(f)$. 

Given a polynomial-like map $f$ with connected Julia set,  
the corresponding {\it polynomial-like germ} is an equivalence class of
polynomial-like maps  $\tl f$ such  that $K(f)=K(\tl f)$ and $f=\tl f$
in a neighborhood of $0$ (hence, by analytic continuation,
also in a neighborhood of $K(f)$). 
We will not make notational distinction between polynomial-like maps and the 
corresponding germs. 
Let
$$
  \mod f = \sup \mod (V\sm U),
$$
where the supremum is taken over all p-l representatives $U\ra V$ of $f$. 

We let $\CC=\CC_d$ be the set of all polynomial-like germs with connected
Julia set.

A unicritical polynomial $p_c:z \mapsto z^d+c$ defines
an element of $\CC$ if and only if $c$ belongs to the Multibrot set
$\MM=\{c \in \C:\, \sup |p_c^n(0)|<\infty\}$.  Those are the only
(normalized) germs with infinite modulus.

We will use superscript $\R$ for the {\it real slice} of a certain space.
For instance $\CC^\R$ stands for germs of
{\it real} polynomial-like maps $f: U\ra U'$ (with connected Julia set),
i.e., such that $f$ preserves the real line and domains
$U$, $U'$ are $\R$-symmetric.

\subsection{Topology}

For a quasidisk $U\subset \C$, let $\BB_U$ stand for the Banach space
of functions holomorphic in $U$ and continuous in $\bar U$. The norm
in this space will be denoted by $\|\cdot \|$.  Let $\BB_U(f,\eps)$
stand for the Banach ball in $\BB_U$ centered at $f$ of radius $\eps$.

We introduce a topology in $\CC$ as follows.
We say that $f_n \to f$ if 
there exists a quasidisk neighborhood $W$ of $K(f)$ such that
(some representatives of) the germs $f_n$ are
defined on $\bar W$ for  sufficiently large  $n$, 
and $f_n$ converges to (an appropriate restriction of a representative of)
$f$ in the Banach space  $\BB_W$. 
The topology in $\CC$ is defined by declaring its closed sets to be
the ones which are sequentially closed.  It is easy to see that
$K(f)$ depends upper semi-continuously
(in the Hausdorff topology) on $f \in \CC$, while its boundary $J(f)$
depends lower semi-continuously.





We let $\CC(\epsilon)$ be the set of all $f\in \CC$ with $\mod(f)
\geq \epsilon$.  Then $\CC(\epsilon)$ is compact, and any compact subset $\KK$ of
$\CC$ is contained in some $\CC(\epsilon)$ (see \cite{McM}).

\subsection{Hybrid classes}
Notice that the Multibrot set $\MM$ has rotational symmetry of order $d-1$ coming from the
fact that polynomials $p_c$ and $p_{\varepsilon c}$ are
affinely equivalent for $\varepsilon= e^{2\pi i/(d-1)}$. In fact,
the {\it moduli space} of unicritical polynomials of degree $d$
(that is, the  space of these polynomial moduli affine conjugacy)
is the orbifold  $\C/ <\varepsilon>$
with order $d-1$ cone point at the origin. 

We say that two
polynomial-like germs $f,g \in \CC$ are {\it hybrid equivalent} if
there exists a quasiconformal map $h:\C \to \C$, such that      
$h \circ f= g\circ h$ in a neighborhood of $K(f)$ (for any representatives
of $f$ and $g$), and $\overline \partial h=0$ on $K(f)$. 
 We call $h$ a {\it hybrid conjugacy} between $f$ and $g$.

By the Douady-Hubbard {\it Straightening Theorem}, 
every $f \in \CC$ is hybrid conjugate to some $p_c$
with $c \in \MM$.  However, in the higher degree case ($d>2$) $c$ may not
be uniquely defined.  Indeed,  polynomials
$p_c$ and $p_{\varepsilon c}$ are affinely equivalent, so they belong
to the same hybrid class.
Vice versa, one can show that hybrid
equivalent polynomials $p_b$ and $p_c$ are
affinely equivalent, so $b=\varepsilon^k c$ for some $k\in \Z/(d-1)\Z$. 
We will see later how to define a single polynomial
straightening associated to each germ (the resolution of the
apparent ambiguity involves global considerations).

We let $\tl \HH_c$ be the hybrid class containing $p_c$.

\subsubsection{Beltrami paths}

A path $f_\lambda \in \CC$, $\lambda \in \D$,
is called a {\it Beltrami path} if there exists a holomorphic
motion
$h_\lambda:\C \to \C$
over $\D$, based on $0$, such that $h_\lambda$ near $K(f_0)$
provides a hybrid conjugacy between $f_0$ and $f_\lambda$.\footnote {The
continuity of $\lambda \mapsto
f_\lambda$ is in fact automatic.}
In this case,
the pair $(f_\lambda,h_\lambda)$ is called a {\it guided Beltrami path}.
The guided Beltrami paths with a fixed initial point $f_0$,
are in one-to-one correspondence with holomorphic families of
Beltrami differentials
$\mu_\la= \dibar h_\la / \di h_\la$ on $\C$ such that  $\mu_0\equiv 0$,
and the differentials $\mu_\la$ vanish a.e. on $K(f_0)$  and are
$f_0$-invariant near $K(f_0)$.
So, in what follows our treatment of Beltrami paths
will freely switch from one point of view to the other.

Obviously any Beltrami path lies entirely in a path connected component of
an hybrid class.

\subsubsection{Hybrid leaves}

Given maps $f_0, f\in \tilde \HH_c$,
let us consider a hybrid conjugacy $h: \C \ra \C$ between them.
Let $\mu$ be the Beltrami differential of $h$
with $L^\infty$-norm $\kappa=\|\mu\|_\infty<1$.  
The family of Beltrami differentials $\la \mu$, $|\lambda|<1/\kappa$, 
generates a guided Beltrami path $(f_\lambda,h_\lambda)$ in
$\tilde \HH_c$, with $f_1$ affinely conjugate to $f$.
In particular each map in $\CC$ can be connected
to one of its straightenings by a Beltrami path.


Let $\HH_c$ be the path connected component of $\tilde \HH_c$ containing
$p_c$.  The $\HH_c$ will be called {\it hybrid leaves}.  By the previous
discussion, $\tilde \HH_c$ is the union of the hybrid leaves
$\HH_{\varepsilon^k c}$, $k \in \Z/(d-1) \Z$.  We will later see that for $c\not=0$, the
hybrid leaves $\HH_{\varepsilon^k c}$, $k \in \Z/(d-1) \Z$ 
(which by definition either coincide or are disjoint), are in
fact all distinct.

\subsection{Expanding circle maps}

A real analytic circle map $g: \T\ra \T$ is called {\it expanding} 
if there exists $n \geq 1$ such that $|Df^n(z)|>1$ for every $z \in \T$. 

Let $\EE=\EE_d$ be the space of real analytic
expanding circle maps $g:  \T \to  \T$ of degree $d$
normalized so that $g(1)=1$.
Such a map admits a holomorphic extension to
a covering $U\ra V$ of degree $d$, 
where $U\ssubset V$ are annuli neighborhoods of $\T$. 
Such extensions will be called {\it annuli representatives of $g$}
and will be denoted by the same letter.  We define
$$
\mod(g)=\sup \mod(V \setminus (U \cup \D))
$$
where the supremum is taken over all annuli representatives $g:U \to V$.

%

Lifting a map  $g\in \EE$ to the universal covering of $\T$, 
we obtain a real analytic function $\tl g: \R\ra  \R$ such that
$\tl g(x)= d\, x + \phi(x)$
where $\phi(x)$ is a $1$-periodic real analytic function with $\phi(0)=0$.
Let $\AA$ be the space of  all such functions, and let
$\AA_n$  be the subspace of the $\phi$  
that  admit a holomorphic extension to the strip $|\Im z| < 1/n$
continuous up to the boundary. 
As the latter spaces are Banach,  $\AA$ is realized
as an inductive limit of Banach spaces,
and we can endow it with  the inductive limit topology.  
It induces a topology on the space $\EE$.
In this topology, a sequence $g_n\in \EE$ converges to $g\in \EE$
if there is a neighborhood $W$ of $\T$ such that all the $g_n$ 
admit a holomorphic extension to $W$, and $g_n\to g$ uniformly on $W$. 

Let $\EE^\R$ stand for the subspace of $\R$-symmetric expanding
circle maps $g: \T\ra \T$ 
(i.e., commuting with the complex conjugacy $z\mapsto \bar z$).   

\begin{lemma} \label {contractible}
The spaces $\EE$ and $\EE^\R$ are contractible.
\end{lemma}

\begin{pf}
  Let us work with the lifts $g: \R\ra \R$ of the maps $g\in \EE$
without making a notational difference between them. 
Let $\EE_1$ stand for the set of  $g\in \EE$ such that $|Dg|>1$ through $\R$. 
This is a convex functional space, so it can be contracted to a point through the affine homotopy. 

The space $\EE_1$ contains the set $\EE_*$ of maps $g\in \EE$ preserving the Lebesgue measure,
so $\EE_*$ can be contracted through $\EE_1$.

To deal with the whole $\EE$,
let us make use of the fact that any $g\in \EE$ has an absolutely continuous invariant measure
$d\mu= \rho\, d\theta$ with real analytic density $\rho(\theta)>0$. 
Let us consider a real analytic circle  diffeomorphism 
$$ 
  h(t)=  \int_0^t \rho(\theta) d\theta 
$$
such that $h_*(d\mu)=d\theta$. Then the map $G= h\circ g\circ h^{-1}$ preserves the Lebesgue measure, so $G\in \EE_*$. 
So, we obtain a projection $\pi: \EE\ra \EE_*$, $g\mapsto G$. 

But the space $\FF$ of diffeomorphisms $h$ is identified with the space of densities $\rho$,
which is also convex, and hence contractible. It follows that $\EE_*$ is a deformation retract for $\EE$,
and the conclusion follows. 

In case of $\EE^\R$,  
just notice that all the above homotopies can be made equivariant 
(with respect to the complex conjugacy). 
\end{pf}



\subsection{External map, mating and product structure}

Given $f \in \CC$, let $\psi:\C \setminus \overline \D \to \C \setminus K(f)$ 
be the Riemann mapping. 
The map $g= \psi^{-1} \circ f \circ \psi$ induces 
(by the Schwarz reflection) an expanding circle endomorphism of degree $d$
called the {\it external map} of $f$.  
It is unique up to conjugacy by a circle rotation, 
so it can be normalized so that $g\in \EE$.  
For $d=2$, this normalization is unique,
but in the higher degree case,
there are generally  $d-1$ ways of normalizing $g$.  Irrespective of this
issue, it is clear that the quality of the analytic extensions of the germ
and of its external map are related by $\mod(f)=\mod(g)$.

\begin{rem}
Maps with symmetries have fewer normalizations, e.g., 
$z\mapsto z^d$ has the maximal possible symmetry group $\Z/(d-1)\Z$ 
and hence has a unique normalization. 
Note that this is the external map of any  polynomial $p_c$, $c\in \MM$. 
\end{rem}

Left inverses to the external map construction are provided by the
{\it matings} between polynomials $p_c$, $c\in\MM$, and expanding
maps $g\in \EE$.  It goes as follows.
Choose a quasiconformal homeomorphism
$h :\C\sm \bar \D \to \C\sm K(p_c)$ such  $h \circ g=p_c\circ h$ near the
circle.
Consider the  Beltrami differential $\mu$ equal to $\dibar h^{-1}/ \di h^{-1}$ on $\C\sm K(p_c)$
and vanishing on $K(p_c)$. 
It is  invariant under $p_c$ on some Jordan disk $D$ containing $K(f)$.
Let $\phi: \C\ra \C$ be the solution of the Beltrami equation
$\dibar \phi/ \di \phi = \mu$.
Then the map  $f = \phi \circ p_c \circ \phi^{-1}$ is polynomial-like on
some neighborhood of $\phi(K(p_c))$, with filled in Julia set $\phi(K(p_c))$, so
up to normalization, it defines a germ in $\CC$.

It is possible to show that, except for the normalization, the germ
$f$ does not depend on the various
choices made in  the construction,
and it clearly depends continuously on $g\in \EE$ and $c\in \MM$.  
In \S \ref {canonic} we will carry out formally
the details of the construction, to obtain the following result:

\begin{thm} \label {productstructure}
There is a canonical choice of
the straightening $\chi(f) \in \MM$ and an external map
$\pi(f) \in \EE$ associated to each germ $f \in \CC$ and
depending continuously on $f$.  It has the following properties: 
\begin{enumerate}
\item For each $c \in \MM$, the
hybrid leaf $\HH_c$ is the fiber $\chi^{-1}(c)$, and the external map $\pi$
restricts to a homeomorphism $\HH_c \to \EE$, whose inverse is denoted by
$i_c$ and called the (canonical) mating,
\item $(\pi,\chi):\CC \to \EE \times \MM$ is a homeomorphism,
\item (Compatibility between matings and Beltrami paths)
For $c,c' \in \MM$, if $f_\lambda$ is a Beltrami path in $\HH_c$ then
$i_{c'} \circ i_c^{-1}(f_\lambda)$ is a Beltrami path in $\HH_{c'}$.
\item External map, straightening and mating are equivariant with
respect to complex conjugation.
\end{enumerate}
\end{thm}

Except for the need to introduce some novelties to handle
$\Z/(d-1)\Z$-ambiguities, the argument follows the quadratic case
\cite {universe}.

One way to understand cancellation of these ambiguities 
(that show up in both the external map and the mating constructions) 
is to introduce {\it markings}. 
Each germ $f \in \CC$ has $d-1$ distinct {\it $\beta$-fixed points}, which
do not disconnect $K(f)$, and a marking of $f$ is just a choice of a
preferred $\beta$-fixed point.  The external map of a marked germ inherits a
marking as well, that is one of its fixed points is distinguished. 
Reciprocally, mating a marked expanding map with a marked
polynomial leads to a well defined marked polynomial-like germ.  Marking also allows us
to resolve the ambiguities inherent to the straightening, since
the straightening of a marked germ is a marked polynomial.

Both expanding maps and polynomials have ``natural markings'',
for expanding maps we choose $1$ as the preferred fixed point, and for
polynomials we choose the landing point of the external ray of angle $0$.
As it turns out, the natural marking of polynomials can be extended
continuously, in a unique way, through the entire $\CC$ 
(this global property is related to the simple topology of $\EE$, see Lemma \ref{contractible}).
Thus, keeping in mind the
natural marking, we end up with natural external map, mating, and
straightening constructions.  By design, the mating construction provides an
inverse to the external map and straightening constructions, so that the
connectedness locus $\CC$ inherits a product structure from $\EE \times
\MM$.

We will discuss markings in more details later in \S \ref {markings} 
(as they are not formally needed for the proof of Theorem \ref {productstructure}).

The reader who is mostly interested in the quadratic case can skip the next two
sections.

\subsection{Proof of Theorem \ref {productstructure}} \label {canonic}

For $c \in \MM$, let $\xi_c:\C \setminus \overline \D \to
\C \setminus K(p_c)$ be the univalent map tangent to the identity at
$\infty$: it satisfies $\xi_c \circ p_0=p_c \circ \xi_c$ on the complement
of $K(p_0)=\overline \D$.

Define the canonical mating $i_c(g)\in \HH_c$ between any polynomial
$p_c$, $c\in\MM$, and any expanding map $g\in \EE$ as follows.
Choose a continuous path $g_t$, $t \in [0,1]$ connecting
$g_0:z \mapsto z^d$ to $g_1=g$, and a continuous family of
quasiconformal maps $h_t:\C \setminus \D \to \C \setminus
\D$, with continuously depending Beltrami differentials $\nu_t$,
satisfying $h_0=\id$ and $h_t \circ g_0=g_t \circ h_t$ near the circle $\T$.
\footnote{Such a family $h_t$ can be constructed as follows.  For large $n$,
$g_{t/n}$ is close to $g_{k/n}$, for every $k \in 0,...,n-1$ and $t \in
[k/n,(k+1)/n]$.  Considering persistent fundamental annuli for the $g_{k/n}$,
one can define a conjugacy $h_{k,t}$, $t \in [k/n,(k+1)/n]$,
between $g_{k/n}$ and $g_{t/n}$
(first on the fundamental domain, then extended by pulling back) satisfying
$h_{k,k/n}=\id$ and the continuity requirements.  Then $h_t$ can be defined
in each interval $[k/n,(k+1)/n]$ by $h_t=h_{k,t} \circ \cdots \circ h_{0,t}$.}
Let $\mu_t$ be the extension of the
Beltrami differential of $h_t \circ \xi_c^{-1}$ to the whole complex plane,
obtained by letting it be $0$ on $K(p_c)$.
It is invariant under $p_c$ in a neighborhood of $K(p_c)$.
Let $\phi_t:\C \to \C$ be the solution of the Beltrami equation
$\dibar \phi_t/\partial \phi_t=\mu_t$.  By invariance of
$\mu_t$, $f_t=\phi_t \circ p_c
\circ \phi_t^{-1}$ is holomorphic in a neighborhood of $K(f_t)=\phi_t(K(p_c))$,
and if $\phi_t$ is appropriately normalized it defines a germ in
$\tilde \HH_c$.  We choose the normalization so that $\phi_t$ depends
continuously on $t$ and $\phi_0=\id$.  

Let us check that  $g_t$ is an external map of $f_t$ for every $t\in [0,1]$.
Let $\psi_t:\C \setminus \overline \D \to \C \setminus K(f_t)$ be the
continuous family of univalent maps normalized so that $\psi_0=\xi_c$, 
and $\tilde g_t :=\psi_t^{-1} \circ f_t \circ \psi_t$
(which extends analytically across the circle by the Schwarz reflection) fixes
$1$, hence $\tilde g_t$ is an external map of $f_t$.
Then $\si_t:= \psi_t^{-1} \circ \phi_t \circ \xi_c:\C \setminus \D \to \C \setminus \D$
is a quasiconformal map conjugating $g_0$ to $\tilde g_t$
whose Beltrami differential coincides with that of $h_t$.  
Hence $\la_t:= \si_t\circ h_t^{-1}$ is a rotation conjugating $g_t$ to $\tl g_t$. 
Since $1$ is a fixed point of $g_t$, $\la_t(1)$ is one of the fixed points of $\tl g_t$ for any $t\in [0,1]$.
But $\la_0=\id$, so  by continuity,
$\la_t(1)=1$ (which is one of the fixed points
of $\tl g_t$) for all $t\in [0,1]$.
We conclude that $\la_t=\id$ and hence  $g_t= \tilde g_t$ for every $t$.

Next, we will  show that the germ $f=f_1$ depends only  on $c$ and $g$, 
but not on the various choices we have made,
which would allow us to define the mating by $i_c(g)=f$.
Let us first show that once the connecting path $g_t$ is chosen,
the path $f_t$ does not depend on the choice of the conjugacies
$h_t$.  Indeed, let $h_t'$ be another choice, with the Beltrami differential $\nu_t'$, 
so that the map $\rho_t:=  h_t^{-1} \circ h'_t$ commutes with $g_0$ near $\T$.
Then the map $\zeta_t  :=\xi_c \circ \rho_t \circ \xi_c^{-1}$ commutes with $p_c$ near $K(p_c)$ (outside it).  
Let us extend $\zeta_t$ to the entire plane by setting $\zeta_t|K(p_c)=\id$.

\begin{lemma} \label {quasiconformality}

The map $\zeta_t$ is a quasiconformal homeomorphism.

\end{lemma}

\begin{pf}

This is a version of the pullback argument, see e.g., \cite{MvS}, Chapter 6,
Section 4.  Choose a quasidisk $V \supset
K(p_c)$ such that $U=p_c^{-1}(V) \ssubset V$ and $\zeta_t \circ p_c=p_c \circ
\zeta_t$ on $U$.  Consider
a continuous family of quasiconformal maps $\zeta^{(0)}_t:\C \to \C$,
such that $\zeta^{(0)}_0=\zeta_0=\id$, $\zeta^{(0)}_t=\zeta_t$ outside $U$,
and $\zeta_t=\id$ near $K(p_c)$.  We can
then set by induction $\zeta^{(k+1)}_t$ as the unique lift (under $p_c$) of
$\zeta^{(k)}_t$ such that
$\zeta^{(k+1)}_t=\id$ near $K(p_c)$.  Clearly $\zeta^{(k)}_0=\id$,
and by continuity in $t$, we see
that $\zeta^{(k)}_t=\zeta_t$ outside $p^{-k}_c(U)$ for every $k$.  
Hence $\zeta^{(k)}_t \to \zeta_t$ pointwise.
Since the dilatations of the $\zeta^{(k)}_t$ do not depend on $k$, 
they form a precompact family of quasiconformal maps.
It follows that the limit map $\zeta_t$ is quasiconformal.
\end{pf}

\begin{rem}

By the Bers Lemma (see \cite {DH}, Lemma 2, p. 303), in order to show that
$\zeta_t$ is quasiconformal, it is enough to check its continuity, i.e., 
that the points $z \in \C \setminus K(p_c)$ near $K(p_c)$
are not moved much by  $\zeta_t$.  
This can be verified directly by an a hyperbolic contraction argument 
(using that $\zeta_t$ commutes with $p_c$): 
in fact, the hyperbolic distance (in the complement of $K(p_c)$) between
$z$ and $\zeta_t(z)$ remains bounded as $z$ approaches $K(p_c)$, see
 (see \cite {DH}, Lemma 1, p. 302).

\end{rem}

We let $\mu_t'=(\xi_c)_*\nu_t'$ outside $K(p_c)$ and $\mu_t'\equiv 0$  on $K(p_c)$. 
Since $\nu_t'= (\rho_t)_*\nu_t$, we have: $\mu_t'= (\zeta_t)_*\mu_t$
(outside  $K(p_c)$ and, obviously, on it). 
Hence $\mu_t'$ is the Beltrami differential for $\phi'_t :=\phi_t \circ \zeta_t$. 
It follows that the map $f_t' :=  \phi'_t \circ p_c \circ (\phi'_t)^{-1}$  
is the mating of $p_c$ and $g_t$ corresponding to the conjugacy $h_t'$. 
But since $\zeta_t$ commuted with $p_c$ near $K(p_c)$, we conclude that $f_t'=f_t$ 
near $K(f_t)$, as was asserted. 
  
Let us now show that the endpoint $f=f_1$
does not depend on the choice of the path $g_t$ connecting $g_0$ and $g$.
Since $\EE$ is simply connected, given another connecting path $g'_t$,
we can fix an homotopy (fixing endpoints) $g^s_t$ with $g^0_t=g_t$ and
$g^1_t=g'_t$.  The mating construction then provides germs $f^s_t \in \CC$,
and it also allows us to choose
hybrid conjugacies $\phi^s_t$ between $p_c$ and $f^s_t$ depending
continuously on $s$ and $t$.
By the previous discussion, $g$ is an external map representative of
$f^s_1$ for every $s$, and in fact
there is a continuous family $\psi^s_1:\C \setminus \overline \D \to \C
\setminus K(f^s_1)$ of univalent maps conjugating $f^s_1$ to $g$. 
Define $\zeta^s:\C \to \C$ by
$\zeta^s=((\psi^s_1)^{-1} \circ \phi^s_1)^{-1} \circ ((\psi^0_1)^{-1} \circ \phi^0_1)$
outside $K(p_c)$ and $\zeta^s=\id$ on $K(p_c)$.
Then $\zeta^s$ commutes with $p_c$ in an outer neighborhood of $K(p_c)$, 
and by the same argument as in Lemma \ref {quasiconformality},
we see that $\zeta^s$ is a global quasiconformal homeomorphism.  
Hence 
$$
 \tau^s:=\phi^s_1 \circ \zeta^s \circ (\phi^s_0)^{-1}=\psi_1^s \circ
(\psi_1^0)^{-1}
$$ 
is a hybrid conjugacy between $f^0_1$ and $f^s_1$ which is also conformal outside $K(f^0_1)$, so it is affine.  
Moreover, note that 1) $\tau^s$ is the identity at $s=0$ and 2) the germs $f^s_1$ are normalized,
so for each $s \in [0,1]$, there is $k \in \Z/(d-1) \Z$ such that
$\tau^s$ is tangent to $z \mapsto e^{2 \pi i k/(d-1)} z$ at $0$. 
We conclude that $\tau^s=\id$  for all $s$.  
Hence $f_1^s=f_1^0$ for all $s$, and in particular $f^1_1=f^0_1$, so the mating $i_c(g)=f$ is indeed well defined.

\begin{lemma}

The mating $(g,c) \mapsto i_c(g)$ is a homeomorphism $\EE \times \MM \to
\CC$.

\end{lemma}

\begin{pf}                                   
Let us begin with continuity of the mating.
It is easy to see that it is continuous in $g$, uniformly with respect to
$c$.
Also, it satisfies $\mod(i_c(g))=\mod(g)$.

So, it is enough to show that for a given $g\in \EE$, 
it is continuous with respect to $c$.
Consider a sequence $c_n \to c\in \MM$.
Choose a path $g_t$ connecting $g_0: z \to z^d$ to $g_1=g$ in $\EE$. 
Passing to a subsequence, we may assume that the
paths $f_{t,n}=i_{c_n}(g_t)$ converge uniformly to a path $f_t$.
 
Then $f_t$ is a path in $\HH_c$. Indeed, since the $\mod(f_{t,n})$ are bounded away from $0$,
the $f_{t,n}$ are $K$-qc conjugate to $f_{0,n}$ (with some $K$ independent of $t$).
Compactness of the space of $K$-qc maps implies that the $f_t$ are $K$-qc conjugate to
$f_0$.  Let us show that $f_t$ is actually hybrid conjugate to $f_0=p_c$.
If this is not the case then $p_c$ must be qc conjugate to
a unicritical polynomial (any straightening of $f_t$)
which is not itself hybrid conjugate to $p_c$, i.e., $c$ is
not qc rigid.  But this implies that $K(p_{c'})$
moves holomorphically for $c'$ in a neighborhood of $c$, $K(p_{c'})= h_{c'} (K(p_c))$.  
It follows that the characteristic function of $K(p_{c_n})$ converges in measure
to that of $K(p_c)$,%
\footnote{To check it, use the following: 
  since $h_c$ as an element of the Sobolev space $W^{1,2}$  depends holomorphically on $c$,
  $\operatorname{Jac}\, h_c = |\di h_c|^2 - |\dibar h_c|^2$ depends continuously  on $c$ weakly in $L^1$. }
 which readily implies that the limit $H$
of hybrid conjugacies $H_n$ between $p_{c_n}$ and $f_{t,n}$ must be a hybrid
conjugacy between $p_c$ and $f_t$ (since $\dibar H_n\to \dibar H$ weakly in $L^2$). 



Let $\psi_{t,n}:\C \setminus \overline \D \to \C
\setminus K(f_{t,n})$ be as in the above construction 
of the external map, i.e., it is
the continuous family of conformal maps such that
$\psi_{t,n}^{-1} \circ f_{t,n} \circ \psi_{t,n}=g_t$ and
$\psi_{0,n}=\xi_{c_n}$.  Then $t \mapsto \psi_{t,n}$ is clearly uniformly
continuous in $t$ (with respect to the topology of uniform convergence on
compact subsets).  So we can take a limiting continuous family $\psi_t:\C
\setminus \overline \D \to \C \setminus K(f_t)$.\footnote {Though
$K(f_{t,n})$ need not converge to $K(f_t)$, any limit is contained in $K(f_t)$ and
its boundary $J(f_{t,n})$ contains $J(f_t)= \di K(f_t)$, which is enough here.}  Then
$\psi_t^{-1} \circ f_t \circ \psi_t=g_t$, and $\psi_0=\xi_c$, 
i.e., $g_t$ is a path that determines the external map of $f_t$.  
So $f_t=i_c(g_t)$ and hence $\lim i_{c_n}(g)=\lim f_{1,n}=f_1=i_c(g)$.  
This proves continuity. 

\msk
Let us now show that the mating is bijective.
Notice first that each polynomial $p_c$ has a single preimage
$(p_0,c)$, since it can only be obtained by mating with its single
external map representative.

Consider a path $f_t$ connecting $p_c$ and an arbitrary
map $f$ in $\HH_c$.  Then $f_t=i_c(g_t)$ where
$g_t$ is the determination of the external map constructed above.
Since $i_c^{-1}(p_c)=\{p_0\}$, this path lifting property implies that
each $i_c:\EE \to \HH_c$ is a bijection.
In particular, $\HH_c=i_c(\EE)$ contains a single polynomial,
$i_c(p_0)=p_c$, so all hybrid leaves are distinct.  Since $\CC$ is the union
of hybrid leaves, this implies that the mating is bijective.

Since the mating is continuous and bijective, it restricts to a
homeomorphism $\EE(\epsilon) \times \MM \to \CC(\epsilon)$ for each
$\epsilon>0$, by compactness, and this implies that the mating is a
homeomorphism $\EE \times \MM \to \CC$.
\end{pf}

We obtain the canonical external map $\pi$ and the canonical straightening
$\chi$ by setting $(\pi,\chi)$ as the inverse of the mating.  All
constructions are clearly equivariant with respect to complex conjugation. 
One also checks directly that
$i_{c'} \circ i_c^{-1}$ takes Beltrami paths in $\HH_c$
to Beltrami paths in $\HH_{c'}$.\qed

\subsection{Marking} \label {markings}

Let us take a polynomial-like map $f: U\ra V$ with connected Julia set.
Let  $A= \bar V \sm U$ and let $\Gamma=\di U$. 
Select an arc $\gamma_0\subset A$ connecting a point $a_{-1}\in f(\Gamma)=\di V$ to one of its preimages, $a_0 \in \Gamma$. 
It can be lifted to an arc $\gamma_1\subset f^{-1} A$ connecting a point $a_1\in
f^{-1}(\Gamma)$ to $a_0$.
In turn, this curve can be lifted to an arc $\gamma_{-2}\subset f^{-2} A $
connecting some point $a_2\in f^{-2} (\Gamma)$ to $a_1$. Continuing this way we obtain a sequence of arcs $\gamma_n$
concatenating a curve $\gamma\in \bar V\sm K(f)$ 
such that $f(\gamma\cap \bar U) =\gamma$ (we will refer to such a curve as
``invariant''). 
A standard hyperbolic contraction argument shows that this curve can
accumulate only at fixed points, and hence lands at some ``preferred''
fixed point of $g$.  Fixed points that arise in this way are called {\it
$\beta$-fixed points}.
A {\it marking} of $g$ is a choice of such an invariant curve up to equivariant
homotopy.  This notion descends to germs, by identifying markings of polynomial-like representatives
of $g$ which coincide up to truncation.

If $c \in \MM$, we can mark the corresponding germ $p_c$
with an {\it invariant external ray} of $p_c$.  No two invariant external rays
can land at the same point,\footnote {Otherwise the sector bounded by
those two rays and which does not contain the critical point would be
invariant by the maximum principle (notice that the image of the sector does
not contain the critical value).} so those markings are indeed distinct.  On
the other hand, it is easy to see that any invariant curve is equivariantly
homotopic to some invariant external ray.  Since such a ray has external argument
$k/(d-1)$ with $k\in \Z/(d-1)\Z$, this procedure shows that there are
exactly $d-1$ different markings, and that markings are in bijection with
the $\beta$-fixed points.

Since an hybrid equivalence between polynomial-like maps gives a correspondence
between the markings, the bijection between markings and $\beta$-fixed points holds
through $\CC$ as well.  The family of $\beta$-fixed points
depends continuously through $\CC$,\footnote {By means of straightening, we
can restrict attention to polynomials, for which it is readily checked that
repelling $\beta$-fixed points are persistent, and repelling
non-$\beta$-fixed points are persistent as well.  Thus a discontinuity might
only arise at a parabolic bifurcation, where both candidates to be a
$\beta$-fixed point are close.} so
each $\beta$-fixed point (or equivalently, each choice of marking) of a germ
admits an unique local continuation to every
sufficient small connected neighborhood a germ.

Similarly to polynomial-like maps,
a circle map $g: U\ra V$ of class $\EE$ can be {\it marked}  by choosing
an invariant curve $\gamma\subset  V\sm \D$ up to equivariant homotopy. 
Such a curve lands at a fixed point of $g$ which depends only on the marking. 
Vice versa, a fixed point determines the marking, 
so there are exactly $d-1$ distinct markings of any circle map $g\in \EE$. 

The marking of $g\in \EE$ corresponding to the fixed point $1$ is called {\it natural}.
It provides us with a continuous global marking of the space $\EE$.  

Due to the product structure $\CC \isom \EE \times \MM$, the natural marking
of $\EE$ can be pulled back to a natural marking of $\CC$.

\subsection{Control of quasiconformal dilatation}

We say that two polynomial-like germs $f,\tl f \in \CC$ are
$(C,\epsilon)$-{\it close} if there exist polynomial-like representatives 
$f:U \to V$ and $\tl f:\tl U \to \tl V$ with 
$\mod(V \setminus  U)>\epsilon$, $\mod(\tl V \setminus  {\tl U})>\epsilon$, 
and a quasiconformal homeomorphism $h:\C \setminus U \to \C \setminus \tl U$
respecting the natural marking of $f$ and $f'$%
\footnote{This condition makes sense since the marking of $f$ can be given by a curve $\gamma\subset \bar V\sm U $
connecting a point $z\in \di U$ to its image $fz\in \di V$ up to homotopy rel the endpoints.}
 with $\Dil(h)<C$ such that
$h \circ f=\tl f \circ h$ on $\partial U$.
Notice that $(C,\epsilon)$-closeness only
depends on the (canonical) external maps $\pi(f)$ and $\pi(\tilde f)$.

Standard arguments (c.f. the proof of Lemma \ref {quasiconformality})
yield:

\begin{lem} \label {closnew}
If $f$ and $\tl f$ are $(C,\eps)$-close by means of $h$
and are hybrid equivalent, then $h$ extends, in a unique way,
to a hybrid conjugacy between $f$ and $\tl f$ with dilatation bounded by $C$.
\end{lem}

By compactness one has:

\begin{lemma} \label {clos}

For every $\epsilon_0>\epsilon>0$ there exists $C>1$
such that if $f,\tl f \in \CC(\epsilon_0)$ then $f$ and $\tl f$ are
$(C,\epsilon)$-close.

\end{lemma}

For nearby germs, the constant $C$ can be taken close to $1$:

\begin{lemma} \label {close}

Let $f_n,\tl f_n \in \CC$ be converging sequences with the same limit.
Then there exists $\epsilon>0$ and $C_n\searrow 1$ such that $f_n$ and $\tl f_n$ are
$(C_n,\epsilon)$-close for every $n$ sufficiently large.

\end{lemma}

\begin{pf}

Let $f=\lim f_n=\lim \tl f_n$.  By definition of convergence, there exists a
polynomial-like representative $f:U \to V$ such that $f_n$ extends
holomorphically to $U$ for every $n$ sufficiently large, $f_n|U$
converges uniformly to $f$, and $f_n \to f$ uniformly on $U$.

Let $W$ be the quasidisk bounded by the equator 
(i.e., the simple closed hyperbolic geodesic) of $V \setminus \overline U$,
and let $\Om:= f^{-1}(W)$,  $\Om_n := f_n^{-1}(W)$ and $\tl\Om_n=\tl f_n^{-1}(W)$.  
Then the Jordan curves $\di \Om_n$ and $\di\tl\Om_n$ converge
in $C^\infty$ topology to the curve $\di\Om$. 
It follows that  for $n$ sufficiently large,
the maps $f_n : \Om_n\ra W$ and $\tl f_n : \tl \Om_n\ra W$ are polynomial-like,
and  the $\mod(W \setminus \Om_n)$, $\mod(W \setminus \tl U_n)$
approach $2\eps:=\mod(W \setminus \Om)$. 
Hence there exist $C^\infty$
diffeomorphisms $h_n:\C \to \C$, such that $h_n|\, \C \setminus W=\id$
and $f_n \circ h_n =\tl f_ n\equiv h_n \circ \tl f_n$ on $\partial \Om_n$, approaching
the identity in the $C^\infty$ topology.  Thus $\Dil h_n\to 1$.  Moreover
$h_n$, being close to the identity, preserves the natural marking and we are
done.
\end{pf}

\subsection{Renormalization and a priori bounds} \label {renorm sec}

A unicritical polynomial-like map $f: U\ra V$ (of degree $d$) is called
{\it renormalizable} with period $p>1$ if there exists a 
topological disk $W\ni 0$ with the following properties:

\ssk\nin R1. The map $g=f^p|\, W$ is a unicritical
polynomial-like map of degree $d$ (onto its image $W'$);
            it is called the {\it pre-renormalization} of $f$.

\ssk\nin R2. The {\it little Julia set} $K(g)$ is connected;

\ssk\nin R3. 
 $K(g)$ does not touch its images
 $f^m(K(g))$, $m=1,\dots, p-1$,
except perhaps at one of its $\beta$-fixed points.   

\ssk
    Note that these images   are also Julia sets  $K(g_m)$ for appropriate  
degree $d$ polynomial-like restrictions $g_m: W_m\ra W_m'$ of $f^p$.
They are also referred to as ``little Julia sets of $f$''.

By \cite[Thm 5.11]{McM}, 
the polynomial-like germ of the renormalization is well defined:
it does not depend on the choice of the domain $W$  above.

In fact, 
there is a standard combinatorial choice of the domain $W$. 
Namely, let us consider the little Julia set $K(g_1)$ around the critical value $f(0)$. 
Among its $\beta$-fixed points, there is a {\it dividing} point $\beta_1$,
i.e., the landing point of more than one external rays (see \cite[Thms 1.2 and 1.4]{Mi}). 
Two of these rays bound a sector containing $f(0)$,
the  {\it characteristic sector} $\SSS_1$. 
The renormalization  range $W_1'$ is obtained by truncating $\SSS_1$ by 
an equipotential and slightly ``thickening'' it
(see \cite{D} or \cite[\S 8]{Mi}).  
The domain $W_1\ni f(0)$ is the pullback of $W_1'$ by $f^p$. 
The domains $W'\supset W\ni 0 $ are the pullbacks of the $W_1' \supset W_1$ under $f$. 

Note that the dividing fixed point $\beta_1$ is uniquely defined.
Indeed, as the characteristic sector $\SSS_1$ has size less than $1/2$,%
\footnote{In fact, $\SSS_1$ is the {\it minimal} sector into which the rays landing on $\orb \beta_1$ divide the plane.}
it does not contain the critical point $0$
and hence $\di \SSS_1$ separates $0$  from  $f(0)$.
  Since the little Julia set $K(g_1)\ni f(0)$ is connected
and the rays landing at the $\beta$-fixed points of $g_1$ do not cut through $K_1(g)$,
there cannot be more than one separating points.

We will mark $\beta_1$ on the little Julia set $K(g_1)$
and the corresponding fixed point $\beta=f^{p-1}(\beta_1)\in f^{-1}(\beta_1)$ on the little Juliua set $K(g)$. 
 Notice that if $f \in \CC^\R$, these points lie in the real line (by symmetry).

\ssk
  Now, the {\it renormalization} of $f$
is obtained by normalizing the pre-renormalization with minimal
possible period,
$$
    Rf(z)= \la^{-1} g(\la z)  : z\mapsto c+z^d+h.o.t.
$$
There is no ambiguity in the choice of normalization 
since the pre-renormalization $g$ is marked with the $\beta$-fixed point.
In case  $f \in \CC^\R$, we have $Rf \in \CC^\R$ as well,
since $\beta$ is real.

%
%
%
%



\msk
The renormalization is called {\it primitive} if the little Julia sets
$K(g_m)$ do not touch, and is called {\it satellite} otherwise.

The set of angles of the external rays
(defined with help of the canonical
straightening of $f$) landing at the distinguished $\beta$-fixed
points of $g$ determine the ``renormalization combinatorics''.\footnote
{An alternative point of view is the following.
The relative positions of the little Julia sets $K(g_m)$ inside the big one,
$K(f)$, can be described in terms of
a graph called the {\it Hubbard tree}.  This graph determines  the
renormalization combinatorics up to symmetry.}
A classical theorem by Douady and Hubbard \cite{DH} asserts that
the renormalizable unicritical  polynomials 
$p_c$ with the same combinatorics form a
``little copy $\MM'$ of the Multibrot set'' (or ``$\MM$-copy''), 
except that the roots of $\MM'$ may or may not be renormalizable.
Thus, the renormalization combinatorics can be labeled by the
little copies themselves.

\ssk
In case of a renormalizable real map,
all the above notions can be described in purely real terms.
The real traces of the little Julia sets are intervals that are permuted
under the dynamics.
The order of these intervals on the line describes the renormalization
combinatorics.
The set of renormalizable maps with a given combinatorics is a parameter
interval $\MM'\cap \R$
called the {\it renormalization window}. Note that the boundary points of a
renormalization
window renormalize to a map with either {\it parabolic}
\footnote{More precisely, it has a parabolic fixed point with multiplier 1.}
or {\it Ulam-Neumann}
\footnote{In this case, $f^2(0)$ is the $\beta$-fixed point
 (such maps are also called {\it Chebyshev}).}
combinatorics. In particular,
{\it the boundary maps are not twice renormalizable}.
(In case of doubling renormalization, the parabolic boundary map is not
renormalizable
in the complex sense, but can be viewed as renormalizable on the real line.)

\ssk A polynomial-like germ $f\in \CC$ is called renormalizable,
if it has a renormalizable representative.
The renormalization descends naturally to the level of germs.  Whether a
germ is renormalizable or not only, and even its renormalization
combinatorics, only depends on its hybrid leaf.  The
renormalization operator acts nicely at the level of hybrid leaves:\footnote
{As for hybrid classes, we notice that affinely conjugated renormalizable
germs have the same renormalization.}

\begin{lemma} \label {paths}

The renormalization operator maps hybrid leaves into hybrid leaves, and takes
Beltrami paths to Beltrami paths.

\end{lemma}

\begin{pf}

It is enough to prove the last statement.  Let $(f_\lambda,h_\lambda)$ be a
guided Beltrami path in a renormalizable hybrid leaf.  Let $f_0:U_0 \to V_0$
be a p-l representative of $f_0$.  We may assume that $V_0$ is small enough
so that $\mu_\lambda=\op h_\lambda/\partial
h_\lambda$ is $f_0$-invariant for every $\lambda \in \D$.
Let $g_0=f_0^p:U'_0 \to
V'_0$ be a pre-renormalization of $f_0$.  Then $\mu_\lambda$ is
$g_0$-invariant for every $\lambda \in \D$.  It follows that $h_\lambda
\circ g_0 \circ h_\lambda^{-1}$ is a pre-renormalization of
$f_\lambda:h_\lambda(U_0) \to h_\lambda(V_0)$.  If $A_0$ is the affine
map conjugating $g_0$ and $R f_0$, there is a unique holomorphic
continuation $A_\lambda$ which normalizes $g_\lambda$, which is readily seen
to conjugate $g_\lambda$ and $R f_\lambda$.  Thus $(R f_\lambda,A_\lambda
\circ h_\lambda \circ A_0^{-1})$ is a guided Beltrami path.
\end{pf}

One can now naturally define $n$ times renormalizable maps,
including $n=\infty$. 
The combinatorics of an infinitely renormalizable map can be labeled  by
a sequence of little Mandelbrot copies $\MM_n'$, $n\in \N$ (describing the combinatorics of the
consecutive renormalizations). It incorporates the sequence $\{p_n\}$ of the 
(relative) renormalization periods. This can be an arbitrary sequence of natural numbers $>1$. 
We say  that $f$ has a {\it bounded combinatorics} if the sequence of periods $p_n$ is bounded.   

We say that an infinitely renormalizable germ $f$ has {\it a priori bounds} if
its renormalizations $R^n f$ have definite moduli: $R^n f\geq \eps >0$.

\comm{
\ssk The bounds are called {\it beau} (over a family $\FF$ of
infinitely renormalizable maps under consideration)
if there exists $\eps_0>0$ such that for any $\de>0$ there exists
moment $n_\de$ such that for any $f\in \FF$ with
$\mod f\geq \de$ we have: $\mod(R^n f)\geq \eps_0$ for $n\geq n_\de$.

The works \cite{K,KL1,KL2} supply a big class of infinitely renormalizable
maps with beau bounds.
In this class the little M-copies $\MM_n'$ describing the combinatorics should stay away from the 
``main molecule''  of $\MM$ (which comprises the main cardioid of $\MM$ and all hyperbolic components obtained from it
via a cascade of bifurcations).  
For instance, this class contains all infinitely renormalizable maps of bounded primitive type and all real infinitely renormalizable
maps with renormalization  periods $p_n\not= 2$.
}

Let us note, for further use, a simple consequence (Lemma \ref {q_n})
of the a priori bounds.  We will need the following
topological preparation:

\begin{lemma} \label {U'}

Let $f':U' \to V'$ be a p-l representative of a
pre-renormalization (not necessarily the first) of
the p-l map $f:U \to V$, of total period $q$.  If $V' \subset V$
then $f^k(U') \subset U$ for $0 \leq k<q$ and $f'=f^q|U'$.

\end{lemma}

\begin{pf}

The connected component
of $f^{-q}(V')$ containing $0$ is a simply connected domain
taken by $f^q$ onto $V'$ as a proper
map which coincides with $f'$ near $K(f')$.  By
analytic continuation, such connected component must coincide with
$U'$ and we have $f'=f^q|U'$.
\end{pf}

\begin{lemma} \label {q_n}

Let $f \in \CC$ be infinitely renormalizable with a priori bounds, and let
$f_n$ be the sequence of pre-renormalizations (of total period $q_n$).
Then there exist $C>0$, $\lambda<1$ (only depending on the a priori bounds)
such that
$$
\max_{m \in \Z/q_n \Z} \diam K_m(f_n) \leq C \lambda^n.
$$

\end{lemma}

\begin{pf}

We are going to show that there exists $\delta>0$, only depending on the
a priori bounds, such that for every
$m,m',n,n'$ such that $n'>n$ and $K_m(f_n) \supset
K_{m'}(f_{n'})$, the Carath\'eodory distance between $K_m(f_n)$ and
$K_{m'}(f_{n'})$ is at least $\delta \diam(K_m(f_n))$.
This clearly implies
that there exists $k>0$ such that
$$
\diam(K_{m'}(f_{n'}))<\diam(K_m(f_n))/2
$$
provided $n' \geq n+k$, and the exponential decay follows.\footnote
{If $\diam(K_{m'}(f_{n'}))>\diam(K_m(f_n))/2$, then we can choose $m_j \in
\Z/\Z_{q_{n+j}}$, $0 \leq j \leq k$ with $m_0=m$ and $m_k=m'$, such that the
$K_{m_j}(f_{n+j})$ are nested.  After suitable translation and
rescaling by $\diam(K_m(f_n))^{-1}$, one gets
$k+1$ compact subsets of the closed unit disk which are pairwise
$\delta/2$-separated.  Since the set of compact subsets of the closed disk
is compact, it is also totally bounded, so $k$ is bounded in
terms of $\delta$.}

Let $f:U \to V$ and $f_{n'}:U' \to V'$ be polynomial-like representatives,
with $\mod(V \setminus U) \geq \mod(V' \setminus U')=\epsilon$.  Up to
replacing $\epsilon$ by $\epsilon/d^t$ and
$U'$ and $V'$ by $f_{n'}^{-t}(U')$ and $f_{n'}^{-t}(V')$ (with $t$ only
depending on $\epsilon$) we may assume that $V' \subset V$.  By Lemma \ref
{U'}, $f'=f^{q_{n'}}|U'$.


Then $f^{q_{n'}}$ has a single critical point in
$\tilde U=f^{m'}(f'^{-1}(U'))$, where we represent
$m'$ in the range $1 \leq m' \leq q_{n'}$.  For each
$1 \leq l \leq q_{n'}$, there exists a unique
$z_l \in K_l(f_{n'})$ such that $f^{q_{n'}-l}(z_l)=0$.  It follows that
$\tilde U$ contains at most one of the $z_l$.  Since $K_m(f_n)$
contains at least two distinct $z_l$, we conclude that $K_m(f_n)
\not\subset \tilde U$.  But
$\mod (\tilde U \setminus K_{m'}(f_{n'})) \geq \epsilon/d$, so $\tilde
U$ is a $\delta \diam(K_{m'}(f_{n'}))$-neighborhood of $K_{m'}(f_{n'})$,
with $\delta$ only depending on $\epsilon$.
\end{pf}

\section{Path holomorphic spaces, the Carath\'eodory metric and the Schwarz
Lemma}\label{Schwarz lemma sec}

A {\it path holomorphic structure} on a space $X$ is a family $\Hol(X)$
of maps $\gamma:\D \to X$, called holomorphic paths, 
which contains the constants and is invariant under holomorphic
reparametrizations: 
for every $\gamma \in \Hol(X)$ and every holomorphic (in the usual sense)
map $\phi:\D \to \D$, $\gamma \circ \phi \in \Hol(X)$.
Natural examples of path holomorphic spaces are
complex Banach manifolds, where holomorphic paths are taken as the paths which
are holomorphic in the usual sense.

If $X,Y$ are path holomorphic spaces, a map
$\Phi:X \to Y$ is called path holomorphic if
for every holomorphic path $\gamma:\D \to X$, $\Phi
\circ \gamma:\D \to Y$ is a holomorphic path.  Let $\Hol(X,Y)$ be the space
of path holomorphic maps from $X$ to $Y$.  Notice that $\Hol(\D,X)=\Hol(X)$.
In case of complex Banach manifolds,
this coincides with the usual notion of being holomorphic,
as long as the maps are continuous.
Obviously, composition of path holomorphic maps is path holomorphic. 


Given a path holomorphic space $X$, any subset $Y \subset X$ can be
naturally considered as a path holomorphic space: if $i:Y \to X$ is the
inclusion, then $\Hol(Y)$ consists of all $\phi: \D\ra Y$ with
$i\circ \phi \in \Hol(X)$.

Let $h(x,y)$ be the hyperbolic metric on $\D$ (normalized to be twice
the Euclidean metric at $0$).
Introduce a metric $d(x,y)$ on $\D$ by taking $d=\frac {e^h-1} {e^h+1}$
(this is a metric by convexity).  
This is the unique metric invariant under the group of conformal automorphisms of $\D$ 
and such that  $d(0,z)=|z|$.

By the usual Schwarz Lemma, any holomorphic map $\phi:\D \to \D$ weakly
contracts $d$: $d(\phi(x),\phi(y)) \leq d(x,y)$.

Let $X$ be a path holomorphic space.
Then we can define the following {\it Carath\'eodory pseudo-metric}:
\be
d_X(x,y)=\sup_{\phi \in \Hol(X,\D)} d(\phi(x),\phi(y)).
\ee
(Obviously, in this definition we can consider only $\phi$  normalized so that $\phi(x)=0$.)
By the usual Schwarz Lemma we have
\be
d_\D(x,y)=d(x,y).
\ee

If $Y \subset X$, we let $\diam_X Y$ denote the diameter of $Y$ in the
pseudo-metric $d_X$. We say that $Y$ is {\it small} in $X$ if $\diam_XY<1$. 

The Carath\'eodory pseudo-metric $d_X$ is a metric if and only if 
bounded path holomorphic functions on $X$ separate points. 
In this case, $X$ is called {\it Carath\'eodory hyperbolic}. 

The definitions immediately imply:

\begin{schw}[weak form]
Any path holomorphic map $\Phi:X \to Y$ is weakly contracting: 
\be
d_Y(\Phi(x),\Phi(y)) \leq d_X(x,y).
\ee
\end{schw}

%

It follows that any subset
of a Carath\'eodory hyperbolic space is Carath\'eodory hyperbolic.

The universal class of
Carath\'eodory hyperbolic spaces are given by Banach balls:  

\begin{lem} \label {emb}
  The unit ball $\BB(1)$ in a complex Banach space $\BB$ is Carath\'eodory hyperbolic and
$d_{\BB_1}(x,0)=\|x\|$.
A path holomorphic space $X$ is Carath\'eodory hyperbolic if and only if $X$ holomorphically
injects into a Banach ball.  
\end{lem}

\begin{pf}
Normalized linear functionals $\phi\in \BB^*_1$ are holomorphic maps $\BB_1\ra \D$.
By definition of the Carath\'eodory metric and the Hahn-Banach Theorem, 
$$
    d_{\BB_1}(x,0) \geq \sup_{\phi\in \BB^*(1)} |\phi(x)| = \|x\|.
$$
The  opposite inequality is obtained by applying the Schwarz Lemma to the embedding
  $\D\ra \BB_1$, $\la\mapsto \la x/\|x\|$ at $\la=\|x\|$. 

The Schwarz Lemma shows that a space which is not Carath\'eodory hyperbolic
can not inject into one that is. 
Vice versa, assume $X$ is Carath\'eodory hyperbolic.  Let $S \subset
\Hol(X,\D)$ be any subset which separates points (e.g., $S=\Hol(X,\D)$).
Then $X$ holomorphically injects into the unit ball of
$\ell^\infty(S)$, the Banach space of bounded functions $S \to \C$,
via the map $x \mapsto (\phi(x))_{\phi \in S}$.
\end{pf}

Small subsets of a hyperbolic space $X$ have  ``definitively stronger''
Carath\'eodory metrics:

\begin{lemma}\label{embeddings are contracting}

Let $X$ be a path holomorphic space and let $Y \subset X$.  Then
for any $x,y\in Y$, 
\be
d_X(x,y) \leq \diam_X(Y) d_Y(x,y).
\ee

\end{lemma}
 
\begin{pf}
Let $r > \diam_XY$. Then for any path holomorphic function $\phi: (X,x)\ra (\D,0)$
we have $\phi(Y) \subset \D_r$.  Hence the function $\tl \phi : =r^{-1} \phi|\, Y$
belongs to $\Hol(Y, \D)$, and we obtain:
$$
    d_Y(x,y)\geq \sup_{\phi} |\tl\phi(y)| = \frac 1r d_X(x,y).
$$
The conclusion follows.
%
\end{pf}

Putting this together with the weak Schwarz Lemma, we obtain:

\begin{schw}[strong form]
Any path holomorphic map $\Phi: Y \to X$ with small image is strongly contracting: 
$$
  d_X(\Phi(x),\Phi(y)) \leq \diam_X(\Phi(Y)) \cdot d_Y(x,y).
$$
\end{schw}

\begin{pf}
  Decompose $\Phi$ as $i\circ \Phi_0$ where $i: \Phi(Y)\ra X$ is the
inclusion, $\Phi_0=\Phi: Y\ra \Phi(Y)$. Then apply the weak
Schwarz Lemma to $\Phi_0$
and Lemma \ref{embeddings are contracting} to $i$. 
\end{pf}

\begin{rem}
  One can consider the (stronger)  Kobayashi metric on path holomorphic spaces as
well.  Though the Kobayashi hyperbolicity is a more general notion,
the spaces of interest in this paper turn out to be already Carath\'eodory
hyperbolic.  More importantly,
the Carath\'eodory metric is much better adapted
to our purposes, since strong contraction can be derived from a very simple
smallness criterion.
\end{rem}

\section{Hybrid leaves as Carath\'eodory hyperbolic spaces} \label {hybrid
classe sec}

\subsection{Path holomorphic structure on hybrid leaves}

For any $c \in \MM$,                                            
we introduce a path holomorphic structure on the hybrid leaf
$\HH_c$ as follows.  
A continuous family $( f_\lambda : U_\la\ra V_\la) \in \HH_c$, $\lambda \in \D$,
is said to be a {\it holomorphic path} if
there exists a holomorphic motion $h_\lambda:\C \to \C$ based at the origin
such that $h_\lambda(K(f_0))=K(f_\lambda))$, 
$\overline \partial h_\lambda=0$ a.e. on $K(f_0)$
(which makes sense since the $h_\la$ are qc by the Quasiconformality
$\la$-Lemma) 
and $h_\lambda \circ f_0=f_\lambda \circ h_\lambda$ on $K(f_0)$ ({\it equivariance} property).

Clearly every Beltrami path is a holomorphic path.  Though the
notion of a Beltrami path is in principle stronger,
they coincide at least locally.  Indeed, let
$f_\lambda$ be a holomorphic path and let $h_\lambda$ be the corresponding
motion of $K(f_\lambda)$.  For each
$\lambda_0 \in \D$, we can make a choice of
a fundamental annulus $V_\lambda \setminus U_\lambda$
which moves holomorphically with $\lambda$ in a small disk $D$ around
$\lambda_0$.  This holomorphic motion can be then extended 
(using the Extension $\la$-Lemma) 
to $\C \setminus U_\lambda$ and then (uniquely) to a 
holomorphic motion on $\C \setminus K(f_\lambda)$ 
that is equivariant on $U_\la\sm K(f)$. 
Matching it with the original motion of $K(f_\lambda)$ 
we obtain a holomorphic motion of $\C$ over $D$,
which provides a hybrid conjugacy.

\begin{rem}

Yet another way to look at holomorphic paths is the following: a continuous
family $f_\lambda \in \HH_c$,
$\lambda \in \D$, is a holomorphic path if and only if
the map $(\lambda,z) \mapsto f_\lambda(z)$ extends to a holomorphic
map in a neighborhood of 
$$
   \cup_{\lambda \in \D} \{\lambda\} \times K(f_\lambda).
$$
However, this point of view
will play no role in our analysis of hybrid classes.

\end{rem}

Through the local characterization of holomorphic paths as Beltrami
paths, we can translate
Theorem \ref {productstructure} (item 3) and Lemma \ref {paths} to path
holomorphicity statements:

\begin{lemma} \label {pat}

\begin{enumerate}
\item All hybrid leaves are path holomorphically equivalent:
For every $c,c' \in \MM$, $i_{c'} \circ i_c^{-1}:\HH_c \to \HH_{c'}$
is path holomorphic.
\item The renormalization operator is leafwise path holomorphic: if $\HH_c$
and $\HH_{c'}$ are such that $R(\HH_c) \subset \HH_{c'}$ then $R:\HH_c \to
\HH_{c'}$ is path holomorphic.
\end{enumerate}

\end{lemma}

\subsection{Carath\'eodory hyperbolicity}

\begin{thm} \label {hyperbolic}

For every $c \in \MM$, $\HH_c$ is Carath\'eodory hyperbolic.

\end{thm}

\begin{pf}

In order to prove Carath\'eodory hyperbolicity of the hybrid leaves, it is
enough, by Lemma \ref {pat}, to prove it for any one of them.  The most
convenient one will be $\HH_0$, since in this case the Julia set traps
a ``definite domain of holomorphicity'':
\be \label {1/4}
\D_{1/4} \subset K(f) \text { if } f \in \HH_0.
\ee
Indeed, there exists a univalent
map $\psi_f$ from $\inter K(f)$ onto $\D$ (the B\"ottcher
coordinate),
such that $\psi_f(f(z))=\psi_f(z)^d$ and
$D\psi_f(0)=1$,\footnote {One way to obtain the B\"ottcher coordinate is
as the restriction of a hybrid conjugacy between $f$ and $p_0$.  The fact
that the derivative at $0$ can be taken as $1$ is immediate from the
normalization.}
and this implies (\ref {1/4}) by the K\"oebe-$1/4$ Theorem.

We will now show that $\HH_0$ holomorphically injects in a Banach ball,
which is equivalent to Carath\'eodory hyperbolicity by Lemma \ref {emb}.  As
target space, we take $\BB_{\D_\rho}$ (the space of bounded holomorphic
functions on $\D_\rho$ which are continuous up to the boundary) for an
arbitrary $0<\rho<1/4$.

Clearly
the restriction operator $I_\rho:\HH_0 \to \BB_{\D_\rho}$
is injective, by analytic continuation.  It is also bounded:
the branch of the $d$-th root of
$f|\inter K(f)$ tangent to the identity at $0$
restricts to a univalent map on $\D_{1/4}$, so that $f|\D_\rho$ can be bounded
in terms of $\rho$.  Let us show that
\be \label {bla3}
f_\lambda(z) \text { is holomorphic in }
(\lambda,z) \in \D \times \D_{1/4} \text { if } f_\lambda \text { is a
holomorphic path in } \HH_0,
\ee
as this clearly implies that $I_\rho$ is path
holomorphic.

Indeed, if $f_\lambda$ is a holomorphic path in $\HH_0$
then there exists a holomorphic motion
$h_\lambda:\C \to \C$ centered on $0$ such that
$h_\lambda(K(f_0))=K(f_\lambda)$, $h_\lambda$ is holomorphic on
$\inter K(f_0)$, and $h_\lambda$ conjugates $f_0$ and
$f_\lambda$
in their filled-in Julia sets.  By separate holomorphicity, we see that
$(\lambda,z) \mapsto (\lambda,h_\lambda(z))$ is holomorphic in
$\D \times \inter K(f_0) \to \bigcup_{\lambda \in \D} \{\lambda\} \times
\inter K(f_\lambda)$.  This implies (\ref {bla3}), since
we can write $f_\lambda(z)=h_\lambda \circ f_0 \circ
h_\lambda^{-1}(z)$.
\end{pf}

\subsection{Carath\'eodory vs Montel}


The following discussion is not actually needed for our proofs of
exponential contraction of renormalization with respect to the
Carath\'eodory metric, but allows us to reinterpret this result in more
familiar terms (the Montel metrics of \cite {universe}).

A compact subset $\KK \subset \CC$ is called {\it sliceable}
if there exist an open quasidisk $W$ and $C>0$ such that $\overline
{\bigcup_{f \in \KK} K(f)} \subset W$ and every $f \in \KK$ has a
holomorphic extension to $W$ bounded by $C$.          
Notice that if
$L \Subset W$ is a  
neighborhood of $0$, then the uniform metric on $C^0(L)$ induces a
distance on $\KK$, and different choices of $L$ lead to
H\"older equivalent distances, by Hadamard's Three Circles Theorem.  In
particular, all those distances define
the same topology on $\KK$, which 
is easily seen to coincide with
the natural topology of $\KK$ (as a subset of $\CC$).

A metric $d$ defined on a compact subset $\KK \subset \CC$
will be called {\it Montel} if
on each sliceable
subset of $\KK$, it is H\"older equivalent to the uniform metric on all
sufficiently small compact neighborhoods of $0$.
Notice that it is enough to check this last condition on any family of
sliceable subsets whose $\KK$-interiors cover $\KK$.  Thus
Montel metrics can be constructed by gluing appropriately metrics on
finitely many sliceable subsets.  They are all H\"older equivalent and
compatible with the topology.

\begin{rem}

We need to go through sliceable subsets
since two different germs $f,\tilde f \in \CC$ may coincide in a neighborhood
of $0$                                                                         
(which cannot happen when  $f$ and $\tilde f$ are in the same sliceable
set).
\end{rem}

\begin{example} Let us sketch  a construction of a pair of
quadratic-like germs $f_\pm \in \CC$ that coincide in a neighborhood of
$0$: indeed $f_+$ and $f_-$ are restrictions of the same
analytic map defined in $K(f_+) \cup K(f_-)$.

Let us start with the map
$$
   F :\T\ra \T,  \quad   z\mapsto \exp (\pi (z-z^{-1})/2),\quad \text{or}\quad 
        \theta\mapsto \pi \sin \theta,
    \quad \text{where}\ z=e^{i\theta}\in \T.
$$
The upper and lower half-circles $\T_\pm $
are invariant under $F$, and the unimodal maps $ F|\, \T_\pm$
admit quadratic-like extensions  that are hybrid equivalent to 
the  Chebyshev map $z \mapsto z^2-2$.

Let us now consider  the analytic real-symmetric  immersion  $\phi:\T \to \C$  satisfying
$$
   \frac {1} {2 \pi^2} \phi(z)^2=1+F(z), \quad \phi(-1)=2 \pi.
$$
Its image  $S$ is a
real-symmetric and $0$-symmetric ``figure eight'' with double point at $0$.
The map $F$ lifts to an analytic  real-symmetric map $f:S \to S $, $f\circ \phi= \phi\circ F$. 
The segments $S_\pm:= \phi(\T_\pm)$ are invarinat under $f$ and 
the maps $f|\, S_\pm$ admit quadraic-like extensions that are hybrid equivalent to 
the Chebyshev map. So, they define two different quadrtic-like germs in $\HH_{-2}$.

However, they define the same germ at $0$ as  $f(z)=2 \pi+z^2+O(z^3)$ as $z\to 0$.


\end{example}

\begin{thm} \label {mon}

For every $\epsilon>0$,
$d_{\HH_c}$ defines a Montel metric on $\HH_c(\epsilon)$,
for each $c \in \MM$.  Moreover those metrics are
uniformly H\"older equivalent to the
restriction of any fixed Montel metric on $\CC(\epsilon)$.

\end{thm}

We will need two preliminary results:

\begin{lemma}

For every $\epsilon>0$, there exists $\delta>0$ such that
$d_{\HH_c}(\delta)$ defines a Montel metric on $\HH_c(\epsilon)$,
for each $c \in \MM$.  Moreover those metrics are
uniformly H\"older equivalent to the
restriction of any fixed Montel metric on $\CC(\epsilon)$.

\end{lemma}

\begin{pf}

By Theorem \ref {hyperbolic} and the Schwarz Lemma,
$\HH_c(\delta)$ is Carath\'eodory hyperbolic, so that $d_{\HH_c(\delta)}$ is
a metric.

Given $\delta>0$, there exists $\rho>0$ such that all $f \in \CC(\delta)$
extends holomorphically to a holomorphic function on $\D_\rho$ bounded by
$\rho^{-1}$.  It then follows, by analytic continuation, that for every $c
\in \MM$ and $z \in \D_\rho$ the function $f \mapsto f(z)$ is
holomorphic on $\HH_c(\delta)$.
This shows that
$2 \rho^{-1} d_{\HH_c(\delta)}$ dominates a Montel metric on
each sliceable subset $\KK \subset \HH_c(\delta)$ (take the Montel metric
given by the uniform distance on each sufficiently small neighborhood $L$
of $0$).

Let us now show that if $\delta>0$ is sufficiently small, then
for each $f \in \HH_c(\epsilon)$,
$d_{\HH_c(\delta)}$ is H\"older dominated by
a Montel metric in an $\HH_c(\epsilon)$-neighborhood
of $f$.  In a neighborhood of $f$ there exist open quasidisks $V
\ssubset V'$ such that if $f_0, f_1 \in \HH_c(\epsilon)$ are
$\gamma$ close with respect to the
Montel metric, then they are $C \gamma^\theta$ close
over $V'$ and there are polynomial-like extensions
$f_0:U_0 \to V$ and
$f_1:U_1 \to V$ with modulus uniformly bounded from below by some $\kappa>0$.
For small $\gamma$, this
implies the existence of a quasiconformal homeomorphism $h:\C \to \C$ which
is the identity outside $V$ and
conjugates $f_0:\partial U_0 \to \partial V$ and $f_1:\partial U_1 \to
\partial V$.  Moreover the Beltrami differential of
$h$ has $L^\infty$ norm bounded by $C' \gamma^{\theta'}$.
If $f_0$ and $f_1$ are hybrid equivalent, we conclude (via the
pullback argument) that $h$ can be turned into a hybrid conjugacy preserving
the natural marking.  This yields a Beltrami path parametrized by
$\D_{C'^{-1} \gamma^{-\theta'}}$ connecting $f_0$ to $f_1$.  This Beltrami
path, restricted to $\D_{C'^{-1} \gamma^{-\theta'}/2}$, lies in
$\HH_c(\kappa/3)$, by the Quasiconformality $\lambda$-Lemma.
Thus,  if $\delta \leq
\kappa/3$ we get $d_{\HH_c(\delta)}(f_0,f_1)
\leq 2 C' \gamma^{\theta'}$, as desired.

The uniformity on $c$ is clear from the argument.
\end{pf}

\begin{lemma}

For every $\epsilon>0$, $d_{\HH_0}$ is a Montel metric on
$\HH_0(\epsilon)$.

\end{lemma}

\begin{pf}

By the Schwarz Lemma $d_{\HH_0}$ is 
dominated by the Montel metric
$d_{\HH_0}(\delta)$ for $\delta>0$ sufficiently small.  For fixed $0<\rho<1/4$,
the restriction $I_\rho:\HH_0 \to \BB_{\D_\rho}$ is a bounded path holomorphic map for
(c.f. the proof of Theorem \ref {hyperbolic}), so another application
of the Schwarz Lemma shows that $d_{\HH_0}$ dominates a multiple of
the uniform metric on $\overline {\D_\rho}$, which is Montel.
\end{pf}

{\it Proof of Theorem \ref {mon}.}  Since both
$d_{\HH_0}$ and $d_{\HH_0(\delta)}$ are Montel
over $\HH_0(\epsilon)$, they are H\"older equivalent.
Since $i_c \circ i_0^{-1}$ is a
biholomorphic map $(\HH_0(\epsilon),\HH_0(\delta),\HH_0) \to
(\HH_c(\epsilon),\HH_c(\delta),\HH_c)$, $d_{\HH_c}$ and
$d_{\HH_c(\delta)}$ are also H\"older equivalent
(with the same constants).
Since the latter is (uniformly) Montel, the former is as well.
\qed

\section{From beau bounds to exponential contraction}\label{beau sec}

{\it A priori bounds} are called {\it beau} (over a family $\FF$ of
infinitely renormalizable maps under consideration)
if there exists $\eps_0>0$ such that for any $\de>0$ there exists
moment $n_\de$ such that for any $f\in \FF$ with
$\mod f\geq \de$ we have: $\mod(R^n f)\geq \eps_0$ for $n\geq n_\de$.

The works \cite{K,KL1,KL2} supply a big class of infinitely renormalizable
maps with beau bounds.
In this class the little $\MM$-copies $\MM_n'$ describing the combinatorics
should stay away from the
``main molecule''  of $\MM$ (which comprises the main cardioid of $\MM$ and
all hyperbolic components obtained from it
via a cascade of bifurcations).
For instance, this class contains all infinitely renormalizable maps of
bounded primitive type and all real infinitely renormalizable
maps with all renormalization  periods $p_n\not= 2$.
Let us emphasize that the approach to the Main Theorem we will
develop in \S\S 6-8 does not rely at all on \cite {K,KL1,KL2} (and it
will cover all real combinatorics, including period doubling, in a unified way).

We will show that beau bounds through
{\it complex} hybrid classes
imply exponential contraction of the renormalization:

\begin{thm} \label {contr}

Let $\FF \subset \CC$ be a family of infinitely renormalizable maps with
beau bounds which is forward invariant under renormalization.  If $\FF$ is
a union of entire hybrid leaves then there exists $\lambda<1$ such that
whenever $f,\tilde f \in \FF$ are in the same hybrid leaf, we have
$$
d_{\HH_{c_n}}(R^n(f),R^n(\tilde f)) \leq C \lambda^n, \quad n \in \N,
$$
where $c_n=\chi(R^n(f))=\chi(R^n(\tilde f))$ and $C>0$ only depends on
$\mod(f)$ and $\mod(\tilde f)$.

\end{thm}

\begin{rem}

We will actually show that
$C(f,\tilde f)$ is small when $f$ is close to $\tilde f$, and indeed
if $\mod(f),\mod(\tilde f) \geq \delta$ we can take
\be \label {C_f,g}
C(f,\tilde f)=C(\delta) d_{\HH_c(\delta)}(f,\tilde f).
\ee

\end{rem}


The proof is based on the
Schwarz Lemma and the following easy ``smallness'' estimate.

\begin{lemma} \label {estimate}

For every $\eps > 0 $ there exist $\de \in (0,\eps)$ and
$\gamma<1$ such that for all $c \in \MM$, we have:
\be
\diam_{\HH_c(\de )} \HH_c(\eps)<\gamma.
\ee

\end{lemma}

\begin{pf}

There exists $r=r(\eps) < 1$ with the following property (see Lemmas \ref
{closnew} and \ref {clos}).
For any p-l germs $f,\tl f \in \HH_c(\epsilon)$, 
there exist p-l representatives $f: U \to V$, $\tl f : \tl U \to \tl V$ and a
hybrid conjugacy (respecting the natural marking)
$h:\C \to \C$ between $f$ and $\tl f$ such that
$\mod(V \setminus U)>\frac {\epsilon} {2}$ and the Beltrami differential $\mu= \dibar h / \di h$
has $L^\infty$-norm bounded by $r(\eps)$. 

Let us consider a Beltrami path  $\D_{\rho} \to \HH_c$,
\be
   \la\mapsto  f_\la=h_{\la \mu} \circ f \circ h^{-1}_{\la \mu}, \quad \mathrm{where}\ 
        \rho= \rho(\eps)=\frac {1+r}{2r}\in (1,\frac 1{r}),
\ee
where $h_{\la\mu}$ is a suitably normalized solution of the Beltrami
equation $\dibar h/ \di h = \la\mu$. 

As $\|\la\mu\|_\infty\leq (1+r)/2$, we have 
$$
  \Dil h_{\la\mu}\leq K=K(\eps)=\frac{r+3}{1-r}, \quad \la\in \D_\rho.
$$
Hence the fundamental annulus of $f_\la$ has modulus at least $\de = \de(\eps):= \eps/2K$,
so $f_{\la} \in \HH_c(\de)$.  By the (weak) Schwarz Lemma,
$d_{\HH_c(\de)}(f,\tl f) \leq d_{\D_\rho} (0, 1) = \rho^{-1}$.
\end{pf}

\noindent{\it Proof of Theorem \ref {contr}.}
Let $\epsilon_0>0$ be the ``beau bound'' for $\FF$, so that for every
$\delta>0$ there exists $n_\delta$ such that $\mod(R^n(f_0)) \geq
\epsilon_0$ whenever $f_0 \in \FF$, $\mod(f_0) \geq \delta$ and $n \geq
n_\delta$.

Using Lemma \ref {estimate}, choose $0<\delta_0<\epsilon_0$ and
$\lambda<1$ such that
$$
\diam_{\HH_c(\delta_0)}
\HH_c(\epsilon_0)<\lambda^{n_{\delta_0}}
$$
for every $c \in \MM$.

The Schwarz Lemma gives for $f, \tilde f \in \HH_c(\delta)$
\begin{align*}
&d_{\HH_{c_n}}(R^n(f),R^n(\tilde f)) \leq \min
\{d_{\HH_c(\delta)}(f,\tilde f),d_{\HH_{c_n}(\delta_0)}(R^n(f),R^n(\tilde
f))\},\\
&d_{\HH_{c_n}(\delta_0)}(R^n(f),R^n(\tilde f)) \leq
d_{\HH_c(\delta)}(f,\tilde f), \quad n \geq
n_\delta,\\
&d_{\HH_{c_n}(\delta_0)}(R^n(f),R^n(\tilde f)) \leq
d_{\HH_{c_{n-n_{\delta_0}}}(\delta_0)}(R^{n-n_{\delta_0}}(f),R^{n-n_{\delta_0}}(\tilde
f)),
\quad n \geq n_{\delta_0}+n_\delta,
\end{align*}
which combined yields the result with $C_{f,\tilde f}$ as in (\ref
{C_f,g}).\qed

\def\Lip{\mathrm{Lip}}
\def\Hol{\mathrm{Hol}}
\def\MI{\mathcal{MI}}

\section{From beau bounds for real maps to uniform contraction}

It is a difficult problem to prove  beau bounds for complex maps.
However, it is more tractable for real maps: in that case,
the beau bounds were established a while ago
(see \cite{S,MvS} for bounded combinatorics and  \cite {LvS,LY} for the general case).

\begin{thm}[Beau bounds for real maps] \label {real beau}

There exists $\epsilon_0 >0$ with the following property.
For every $\delta>0$ there exist
$\eps=\eps(\de)>0$ and $N=N(\delta)$ such that for any  $f \in \CC^\R(\delta)$,
we have 
$$
  R^n(f) \in \CC(\epsilon), \ n=0,1,\dots \quad \mbox{and}\quad R^n(f) \in \CC(\epsilon_0) \   n= N, N+1,\dots .
$$

\end{thm}

From the point of view of the previous discussion, the main shortcoming of
this result is that it does not provide enough compactness for complex maps,
which is crucial for the Schwarz Lemma application.
In this section we will  overcome this problem
using some ideas from functional analysis and differential topology 
(Theorem \ref {real beau} remaining the only ingredient from ``hard
analysis''), by proving:

\begin{thm}[Beau bounds and macroscopic contraction for complex maps in the real hybrid classes]\label{beau complex}
There exists $\eps_0>0$ with the following property. 
For any $\gamma >0$ and $\de>0$ there exists $N=N(\gamma,\de)$ such that for any two maps 
$f, \tl f \in \CC(\de)$ in the same real-symmetric hybrid leaf we have
$$
        R^n f, R^n \tl f\in \CC(\eps_0), \quad {\mathrm and}\quad
d_{\HH_{c_n}}(R^n f, R^n \tl f) < \gamma, \quad n\geq N,
$$
with $c_n=\chi(R^n f)=\chi(R^n \tilde f)$.
\end{thm}

The proof of this result will take this and the next two sections.

Through the sequel,
let $\II$ be the subspace of infinitely renormalizable p-l maps in $\CC$, and let 
$\II(\de)= \II\cap \CC(\de)$.
Let $\CC^{(\R)}$
 be the space of polynomial like maps with connected Julia set
which  are hybrid equivalent to real p-l maps.
We will use superscript $(\R)$ for the slices of various spaces by $\CC^{(\R)}$,
e.g.,  $\II^{(\R)}(\de):=\II\cap \CC^{(\R)}(\de) $.

%
%

\subsection{Cocycle setting}

We will now abstract properties of the renormalization operator 
that will be sufficient for Theorem \ref{beau complex}.

\def\SSS{{\mathcal{S}}}

Let $\SSS$ be a semigroup, and let $Q=\{(m,n)\in
\N\times \N:\ n>m\}$.
An $\SSS$-{\it cocycle} is a map $G:  Q \ra \SSS$, $(m,n)\mapsto G^{m,n}$,
such that
$G^{m,n} G^{l,m}  = G^{l,n}$.

Letting $F_n:= G^{n, n+1} \in \SSS$, we obtain
\begin{equation}\label{cocycle}
    G^{m,n}= F_{n-1} \circ \cdots \circ F_m,
\end{equation}
and vice versa, any sequence $F_n\in \SSS$ determines a cocycle by means of
(\ref{cocycle}).


Let $\Hol(\HH_0,\HH_0)$ be the semigroup of continuous path holomorphic maps
$F:\HH_0 \to \HH_0$.
Let $\Hol^\R(\HH_0,\HH_0)$ be the sub-semigroup of those $F$
such that $F(\HH^\R_0) \subset \HH^\R_0$.

\begin{thm} \label {ff}

Let $\GG$ be a family of cocycles with values in
$\Hol^\R(\HH_0,\HH_0)$ satisfying:
\begin{enumerate}


\item [H1.] {\it A priori bounds for complex maps}:  
For every $\de> 0$, there exists $\eps=\eps(\de) >0$ such that
$$
\mbox{If $f \in \HH_0(\de)$  
then  $G^{m,n}(f) \in \HH_0(\epsilon)$ for every $G \in \GG$ and $(m,n)\in Q$.}
$$

\item [H2.] 
{\it Beau bounds for nearly real maps}:
There exists $\epsilon_0>0$ such that for every $\delta>0$, there
exists $N=N(\delta)$ and $\eta=\eta(\delta)>0$ such that: 
$$
\mbox{If  $f \in \HH_0^\R(\delta)$, $\tl f \in \HH_0(\delta)$ and $d_{\HH_0}(f,\tl f) \leq \eta$}
$$
$$
\mbox{ then $G^{m,n}(\tl f) \in \HH_0(\epsilon_0)$ for every $G \in \GG$ and
$(m,n) \in Q$ with $n-m\geq N$.}
$$
\end{enumerate}
Then we have:
\begin{enumerate}
\item [C1.] {\it Macroscopic contraction}: For every $\de>0$ and $\gamma>0$ there exists $N=N(\de,\gamma)$
such that if $f,\tl f \in \HH_0(\de)$ then
$$
\mbox{ $d_{\HH_0}(G^{m,n}(f),G^{m,n}(\tl f)) < \gamma$ for every
$G \in \GG$ and $(m,n) \in Q$ with $n-m \geq N$.} 
$$

\item [C2.] {\it Beau bounds for complex maps}: 
 For every $\delta>0$ there exists
$N=N(\delta)$ such that
$$
\mbox{ $G^{m,n}(\HH_0(\delta)) \subset \HH_0(\epsilon_0)$ for every $G \in
\GG$ and $(m,n) \in Q$ with $n-m\geq N$.}
$$
\end{enumerate}

\end{thm}

\subsection{Reduction to the cocycle setting}




Let $\Pi=i_0 \circ \pi:\CC \to \HH_0$, so that for each $c \in \MM$, $\Pi$
restricts to a homeomorphism $\HH_c \to \HH_0$ which is path holomorphic and
preserves the modulus.
For each infinitely renormalizable hybrid leaf $\HH_c$,
we can associate a cocycle $G=G_c$ with values in $\Hol(\HH_0,\HH_0)$, by
the formula
\be
G^{m,n}(\Pi(f))=\Pi(R^{n-m}(R^m(f))), \quad f \in \HH_c.
\ee

Let $\GG$ be the set of all such cocycles which correspond to
real-symmetric hybrid leaves.
Once we show that hypothesis (H1-H2) of Theorem \ref {ff} are satisfied for
$\GG$, the conclusions  (C1-C2) translate precisely into beau bounds and
macroscopic contraction in real hybrid classes (Theorem~\ref {beau complex}).

Let us start with (H1): 

\begin{lemma}

For every $\de>0$ there exists 
$\epsilon=\eps(\de) >0$ such that:
$$\mbox{ If $f \in \II^{(\R)}(\de)$ then
$R^n f \in \CC(\epsilon)$, $n=0,1,\dots $.}
$$
\end{lemma}

\begin{pf}

Let $g = p_c: z\mapsto z^d+c$, $c\in \R$,  be the straightening of $f$,
 and let $g_n$ denote its $n$th pre-renormalizations.
Let $q_n$ be the corresponding periods, that is, $g_n$ is a p-l restriction of $g^{q_n}$.

By Lemmas \ref {closnew} and \ref {clos},
there exist p-l representatives $g:U' \to V'$, $f:U \to V$ with
$\mod(V' \setminus U'),\mod(V \setminus U)>\de/2$ and a $C$-qc map
$h:(\C, U) \to (\C, U')$ with
$C=C(\de)$ conjugating $f$ to $g$, 
$h \circ f= g \circ h$ in $U$. 
Since $\mod(V' \setminus U')>\delta/2$,
$$
\inf_{y \in \partial U'} \inf_{x
\in K(g)} |y-x| \geq A \diam K(g),
$$
for some $A=A(\delta)>0$.

By the {\it a priori  bounds} for real maps, 
there exists $\eta>0$, depending only on the degree $d$, such that
$$
  \mod g_n \geq \eta, \quad n=0,1,\dots.  
$$
It follows that each germ $g_n$ has a  (non-normalized)
p-l representative $U_n \to V_n$ with 
$\mod( V_n \setminus U_n) >\eta'$ and 
$$
\sup_{y \in \partial g^k(U_n)} \inf_{x \in g^k(K(g_n))} |y-x|
\leq A \diam g^k(K(g_n))\leq A \diam K(g),
\quad  k =0,1,2,\dots, q_n-1
$$
with $\eta'=\eta'(A,\eta)>0$, by compactness of $\CC(\eta)$.

\comm{
By Lemmas \ref {closnew} and \ref {clos},
there exist p-l representatives $g:U' \to V'$, $f:U \to V$ with
$\mod(V' \setminus U'),\mod(V \setminus U)>\de/2$ and a $C$-qc map
$h:(\C, U) \to (\C, U')$ with
$C=C(\de)$ conjugating $f$ to $g$, 
$h \circ f= g \circ h$ in $U$. 
Since $\mod(V' \setminus U')>\delta/2$,
$$
\sup_{y \in \partial U'} \inf_{x
\in K(g)} |y-x| \leq A \diam K(g),
$$
for some $A=A(\delta)>0$.
}

We conclude that $g^k(U_n) \subset U'$ for $0 \leq k \leq q_n-1$.
Consequently,  the p-l map 
$$ 
  h^{-1} \circ g_n \circ h:\ h^{-1}(U_n) \to h^{-1}(V_n)
$$ 
is a (non-normalized) 
representative of the $n$-th pre-renormalization of $f$.  Since
$$ \mod(h^{-1}(V_n) \setminus h^{-1}(U_n))>\eta' / C, $$ 
we obtain {\it a priori} bounds for $f$  with $\epsilon = \eta' /
C$.  
\end{pf}

In order to show that hypothesis (H2) is satisfied, we will use the following.

\begin{lemma}

Let $f_n, \tl f_n \in \II$ and $\chi(f_n)=\chi(\tl f_n)$. 
Assume that the sequences $f_n$ and $\tl f_n$  converge to the same limit $f$. 
If $k_n \to \infty$ then
$\liminf_{n \to \infty} \mod(R^{k_n}(f_n))=\liminf_{n \to \infty} \mod(R^{k_n}(\tl f_n))$.

\end{lemma}

\begin{pf}

It is enough to show that for every $\epsilon>0$,
$$
\text {if} \quad
\liminf_{n \to \infty} \mod(R^{k_n}(f_n))>\epsilon \quad \text {then}
\quad \liminf_{n \to \infty} \mod(R^{k_n}(\tl f_n)>\epsilon.
$$

Let $f_n: U_n \ra V_n$ and $\tl f_n : \tl U_n\ra \tl V_n$ 
be p-l representatives of the germs $f_n$ and $\tl f_n$ that  Carath\'eodory converge to
a p-l map $f: U\ra V$. 
By Lemma \ref {close},
there exist $C_n\to 1$ and  $C_n$-qc maps $h_n : V_n\ra \tl V_n$ conjugating $f_n$ to $\tl f_n$
(maybe after a slight adjustment of the domains). 

Let $f'_n: U_n' \to V_n'$ be p-l representatives of the
$k_n$-th pre-renormalizations of the $f_n$
with $\liminf \mod(V_n' \setminus U_n')>\epsilon$
and filled Julia sets $K_n$.  By Lemma \ref {q_n}, $\diam(K_n) \to 0$, so
we may choose $V_n'$ contained in $V_n$.  By Lemma \ref {U'},
$f'_n=f_n^{q_{k_n}}|U'_n$.

Let $\tl U_n'=h_n(U_n')$ and $\tl V_n'= h_n(V_n')$.
Then the map $\tl f_n^{q_{k_n}}: \tl U_n'\ra \tl V_n'$ is a well defined  
p-l representative of the $k_n$-th pre-renormalization of $\tl f_n$.   
Moreover, $\mod (\tl V_n'\sm \tl U_n')>\mod(V_n' \setminus U_n')/C_n$,
and the conclusion follows. 
\end{pf}

Take  $\epsilon_0$  from Theorem \ref {real beau}.
If (H2) does not hold for $\eps_0/2$, we can find a $\delta>0$ and
sequences $f_n \in \II \cap \CC(\delta)$, $\tl f_n \in \HH_{\chi(f_n)}(\delta)$, 
$k_n \to \infty$, with $d_{\HH_{\chi(f_n)}}(f_n, \tl f_n) \to 0$
such that $R^{k_n} \tl f_n \notin \CC(\epsilon_0/2)$.
Passing to a subsequence, we can assume that the $f_n$ converge, 
and thus the $\tl f_n$ must converge to the same limit.
By the previous lemma,
it follows that $\liminf \mod(R^{k_n} f_n) \leq \epsilon_0/2$,
contradicting Theorem~\ref {real beau}.

We have reduced Theorem \ref {beau complex} to Theorem \ref {ff}.

\subsection{Retractions and the proof of Theorem \ref {ff}}

Recall that a continuous map $P: X\ra X$ in a topological space $X$ is
called a {\it retraction} if $P^2=P$.
In other words, there exists a closed subset $Y\subset X$ (a ``retract'')
such that
$P(X)\subset Y$ and $P|\, Y =\id$. Linear retractions in topological vector
spaces are
called projections.  A retraction is naturally called
{\it constant} if its image is a single point.

The proof of Theorem \ref {ff} has two main parts.  The first shows,
using (H1), that lack of uniform contraction in the leafwise dynamics
allows one to construct a retraction towards a non-trivial
``attractor'':

\begin{thm} \label {trivial}

Let $\GG$ be a family of cocycles with values in $\Hol(\HH_0,\HH_0)$.
Assume that property (H1) holds but (C1) fails.  Then there exist sequences
$G_k \in \GG$, $(m_k,n_k) \in Q$ with $n_k-m_k \to \infty$, and a
non-constant retraction $P \in \Hol(\HH_0,\HH_0)$ such that
$G_k^{m_k,n_k}(f) \to P(f)$ for every $f \in \HH_0$.

\end{thm}

The second shows, in general, that non-trivial retractions
can not be ``too compact'' in the real direction.

\begin{thm} \label {retraction}

Let $P \in \Hol^\R(\HH_0,\HH_0)$, and assume the following
compactness property for nearly real maps:

\ssk
(P) {\it There exists a compact set $\KK \subset \HH_0$ such that if
$f_n \in \HH_0$ is a sequence converging to $f \in \HH_0^\R$, then
$P(f_n) \in \KK$ for $n$ large.}
\ssk

Then $P$ is constant.

\end{thm}

Those two results put together are then seen to imply Theorem \ref {ff}:

\noindent{\it Proof of Theorem \ref {ff}.}
Let $\GG$ be a family of cocycles with values in $\Hol^\R(\HH_0,\HH_0)$
such that
(H1) and (H2) hold, but (C1) does not.
By Theorem \ref {trivial}, there exist a
sequence $G_k \in \GG$ and $(m_k,n_k)\in Q$ with $n_k-m_k \to \infty$, and a
non-constant retraction $P \in \Hol(\HH_0,\HH_0)$, such that
$G^{m_k,n_k}_k(f) \to P(f)$ for every $f \in \HH_0$.
It satisfies the hypothesis of Theorem \ref {retraction}:
since each $G^{m_k,n_k}_k$ preserves $\HH_0^\R$, $P$ also does, while (H2)
immediately gives\\
\ssk
 (P') {\it For any $f \in \HH^\R_0(\de)$ and $\tl f \in \HH_0(\de)$ with
$d_{\HH_0}(\tl f,f)<\eta$
         we have $ P(\tl f) \in \HH_0(\epsilon_0)$}\\
\ssk
which clearly implies the crucial property (P).  By Theorem \ref
{retraction}, $P$ is
constant, yielding the desired contradiction.

\ssk
This concludes the proof of (C1).  Together with (H2), it implies (C2).
\qed

The essentially independent
proofs of Theorems \ref {retraction} (of differential topology nature)
and \ref {trivial} (dynamical) will be given in the
next two sections.

\section{Triviality of retractions}\label{triviality sec}

\subsection{Plan of the proof of Theorem \ref {retraction}} \label {planofproof}

Let us describe the plan of the proof of Theorem \ref {retraction}.
Let $P^\R= P|\, \HH_0^\R$.
Let  $\ZZ^\R = \mathrm {Im}\, P^\R = \Fix P^\R$.
By Property (P),  $\ZZ^\R \subset \KK$,
and hence $\ZZ^\R = P(\KK)$, which  is compact.
We will complete the argument in three consecutive steps:

\begin{itemize}
\item [{\rm Step 1.}]  $\ZZ^\R$ is a finite-dimensional topological manifold
(by the Implicit Function Theorem),
\item [{\rm Step 2} ] $\ZZ^\R$  is  a single point
(by a Brower-like topological argument),
\item [{\rm Step 3.}] $\ZZ:=\mathrm {Im} P=\Fix P$ is a single point, too (by
analytic continuation).
\end{itemize}

\comm{
(P) {\it There exists $\epsilon_0>0$ such that for every $\delta>0$ there
exists $\eta=\eta(\delta)>0$ such that if
$f \in \HH^\R_0(\de)$ and $\tl f \in \HH_0(\de)$ are such that
$\dist_{\HH_0}(\tl f,f)<\eta$ then $P(\tl f) \in \HH_0(\epsilon_0)$.}
\ssk
}

The first and third steps would be immediate to carry out if we were dealing
with Banach spaces.  For instance, corresponding to the first step we have:

\begin{lemma}[\cite{Ca}] \label {dp}

Let $\BB$ be a complex (respectively, real) Banach space and let
$P$ be a holomorphic (respectively, real analytic)
map from an open set of $\BB$ to $\BB$ such
that $P(0)=0$.  Assume that $DP(0)$ is compact and
$P^2=P$ near $0$.  Then for any sufficiently
small open ball $B$ around $0$ in
$\BB$, $P(B)$ is a complex (respectively, real analytic)
finite-dimensional manifold.

\end{lemma}

\begin{pf}

This is a particular case of \cite {Ca} but we will give the argument for
the
convenience of the reader.  Let $h=\id-DP(0)-P$.  Since $P^2=P$, we
have $DP(0)=DP(0)^2$ and hence $Dh(0)^2=\id$, so $h$ is a local
diffeomorphism
near $0$.  Obviously $h \circ P=DP(0) \circ h$, so $P(B)=h^{-1}(DP(0)(h(B)))$.
Since $DP(0)$ is compact and $DP(0)=DP(0)^2$, it
has finite rank so $DP(0) \cdot h(B)$ is an open subset of a finite dimensional
subspace.
\end{pf}

In order to translate this more familiar analysis to
our context, we will use {\it Banach slices} (first introduced in \cite
{universe}).

\subsection{Banach slices} \label {bs}

Let $f \in \HH_0$.  We call an open quasidisk $W$ $f$-admissible if
$W \supset K(f)$ and $f$ extends holomorphically to $W$ and continuously to
the boundary.  Let $\BB^*_W \subset \BB_W$
be the Banach space
of all holomorphic maps $w:W \to \C$ such that $w(z)=O(z^{d+1})$
near $0$ and $w$ extends continuously to the boundary.

The set of $f$-admissible quasidisks can be partially ordered
by inclusion.  This partial order is directed in the sense that
any finite set has a lower bound,
thus it makes sense to speak of ``sufficiently small'' $f$-admissible $W$.

Let $B_{W,r}$ be the
open ball around $0$ in $\BB^*_W$ of radius $r$.  The following lemma is a
straightforward consequence of the definition of the topology in $\HH_0$.

\begin{lemma}

For every $\epsilon>0$, and $f \in \HH_0(\epsilon)$, for every sufficiently
small $f$-admissible $W$, for every $r>0$, there exists a neighborhood
$\VV$ of $f$ in $\HH_0(\epsilon)$ such that for every $f' \in \VV$,
$W$ is $f'$-admissible and $f'-f \in B_{W,r}$.

\end{lemma}

It is easy to see that there exists $\epsilon_0=\epsilon_0(f,W)$ and
$r_0=r_0(f,W)>0$ such that if $w \in B_{W,r_0}$ then
$f'=f+w$ admits a polynomial-like restriction $f':U \to V$ with $K(f')
\subset U \subset W$ and $\mod(V \setminus U)>\epsilon_0$.  Since $f'(0)=0$,
$f'$ defines a germ in $\HH_0(\epsilon_0)$ denoted by $j_{f,W,r_0}(w)$.
The map $j_{f,W,r_0}:B_{W,r_0} \to \HH_0$ is readily seen to be continuous
and injective.

\begin{thm} \label {diff}

Let $f$, $W$ and $r_0$ be as above, and let
$\lambda \mapsto w_\lambda$ be a continuous map $\D \to B_{W,r_0}$.
Then $w_\lambda$ is holomorphic if and only if
$f_\lambda=j_{f,W,r_0}(w_\lambda)$ is a holomorphic path in $\HH_0$.

\end{thm}

For the proof, we will need a preliminary result.  As in the proof of
Theorem \ref {hyperbolic}, for $0<R<1/4$ we let
$I_R:\HH_0 \to \BB_{\D_R}$ be the restriction operator, which is well
defined by (\ref {1/4}).

\begin{lemma} \label {crit}

Let $0<R<1/4$ and
let $\lambda \mapsto f_\lambda$ be a continuous map $\D \to \HH_0$.  Then
$f_\lambda$ is a holomorphic path in $\HH_0$
if and only if $\lambda \mapsto I_R(f_\lambda)$ is a holomorphic path in
$\BB_{\D_R}$.

\end{lemma}

\begin{pf}

The only if part (equivalent to the path holomorphicity of $I_R$) was
established in the proof of Theorem \ref {hyperbolic}.

Assume that $I_R(f_\lambda)$ is a holomorphic path in $\BB_{\D_R}$.
Since $f_\lambda$ is assumed to be continuous, in order to show that
it is a holomorphic path, it is enough to
construct a holomorphic motion $h_\lambda:\inter K(f_0) \to \inter
K(f_\lambda)$ such that for each $\lambda$, $h_\lambda$ is holomorphic and
conjugates $f_0$ and $f_\lambda$: by the
Extension $\lambda$-Lemma Theorem,
it extends to a holomorphic motion $\C \to \C$, which, by
continuity, conjugates
$f_0|K(f_0)$ and $f_\lambda|K(f_\lambda)$.
For the construction, we will make use of the B\"ottcher coordinate (c.f.
proof of Theorem \ref {hyperbolic}) $\psi_f:\inter K(f) \to \D$ associated
to any map $f \in \HH_0$:
the desired holomorphic motion is then given by
$h_\lambda=\psi_{f_\lambda}^{-1} \circ \psi_{f_0}$.  It is
obviously injective
and holomorphic in $z$, and conjugates $f_0$ and $f_\lambda$, for each
$\lambda \in \D$, so we just need to
show that $\psi_{f_\lambda}^{-1}(z)$ is a holomorphic function $\D \times \D
\to \C$.

Holomorphicity of $I_R(f_\lambda)$ implies that
$f_\lambda(z)$ is a holomorphic function $\D \times \D_R \to \C$.
By the K\"oebe-$1/4$ Theorem, if $z \in \D_{R^3}$ then
$\psi_{f_\lambda}(z) \in \D_{R^2}$, so $\psi_{f_\lambda}(z)^d \in
\D_{R^4}$, and (again by the K\"oebe-$1/4$ Theorem),
$f_\lambda(z)=\psi_{f_\lambda}^{-1}(\psi_{f_\lambda}(z)^d) \in \D_{R^3}$.
It follows that $f_\lambda^n(z) \in \D_{R^3} \subset
\D_R$ for every $z \in \D_{R^3}$, $n \geq 1$.
By holomorphic iteration we conclude that for every $n \geq 1$,
$f_\lambda^n(z)$ is holomorphic in $(\lambda,z) \in \D \times
\D_{R^3}$.

Let $\psi_{f_\lambda,n}:\inter K(f_\lambda) \to \C$ be such that
$\psi_{f_\lambda,n}(z)^{d^n}=f_\lambda^n(z)$ and $D
\psi_{f_\lambda,n}(0)=1$.  It is easy to see that $\psi_{f_\lambda,n}$
converges to $\psi_{f_\lambda}$ uniformly on compacts of $\inter
K(f_\lambda)$ (actually one usually constructs the
B\"ottcher coordinate $\psi_{f_\lambda}$
directly as the limit of the $\psi_{f_\lambda,n}$).  Over $(\lambda,z) \in
\D \times \D_{R^3}$, the holomorphicity of $f_\lambda^n(z)$
implies, successively, that $\psi_{f_\lambda,n}(z)$ and
$\psi_{f_\lambda}(z)$ are also holomorphic.

By the K\"oebe-$1/4$
Theorem, $\psi_{f_\lambda}^{-1}(\D_{R^4}) \subset \D_{R^3}$, and it
follows that $\psi_{f_\lambda}^{-1}(z)$ is a holomorphic function of
$(\lambda,z) \in \D \times \D_{R^4}$.  Since for each fixed $\lambda \in
\D$, $\psi_{f_\lambda}^{-1}$ is a univalent function of $\D$,
Hartog's Theorem implies that $\psi_{f_\lambda}^{-1}(z)$ is in fact
a holomorphic function of $(\lambda,z)$ through $\D \times \D$.
\end{pf}

\noindent{\it Proof of Theorem \ref {diff}.}
Assume that
$f_\lambda$ is a holomorphic path in $\HH_0$,
and let us
show that for every bounded linear functional $L:\BB^*_W \to \C$, $\lambda
\mapsto L(w_\lambda)$ is holomorphic: since $w_\lambda$ takes values in a
ball, this implies that
$\lambda \mapsto w_\lambda$ is holomorphic.  By (\ref {1/4}), $f_\lambda(z)$,
and hence $w_\lambda(z)=f_\lambda(z)-f(z)$,
is holomorphic in
$\D \times \D_{1/4}$.  By Hartog's Theorem, it is then
holomorphic in $\D \times W$, and since it is continuous
in $z$ up to $\partial W$, and bounded on both variables,
we see that $\lambda \mapsto f_\lambda(z)$ is a
holomorphic function for every $z \in \overline W$.  By Riesz's Theorem,
there exists a complex measure of finite mass $\mu$,
supported on $\overline W$, such that $L(w)=\int w(z) d\mu(z)$, so
$\lambda \mapsto L(w_\lambda)$ is holomorphic.

Assume now that $\lambda \mapsto w_\lambda$ is holomorphic.  Since
$j_{f,W,r_0}$ is continuous, $\lambda \mapsto f_\lambda$ is continuous as
well.  Fix $0<R<1/4$.
By (\ref {1/4}), $\D_{1/4} \subset K(f) \subset W$, hence the restriction
operator $I_{R,W}:\BB^*_W \to \BB_{\D_R}$ is holomorphic. 
Since $w_\lambda \in \BB^*_W$ depends
holomorphically on $\lambda$, it follows that
$I_R(f_\lambda)=I_R(f)+I_{R,W}(w_\lambda) \in \BB_R$ also
depends holomorphically on $\lambda$.
By Lemma \ref {crit}, $f_\lambda$ is a holomorphic path.
\qed

\subsection{Proof of Theorem \ref {retraction}}

We will carry out the three steps of the
plan of proof described in \S \ref {planofproof}.  We will use the notation
introduced therein.

Let $\epsilon>0$ be such that $\KK \subset \HH_0(\epsilon)$.
Let $f \in \ZZ^\R$.
Let us consider a  neighborhood $\UU$ of $f$ in $\HH_0(\eps/2)$.
If it is small enough then all the maps $f\in \bar \UU$
are well defined on some admissible neighborhood of $W\supset K(f)$ (see \S
\ref {bs}),
so $\bar\UU$ naturally embeds into some Banach ball $B_r:= B_{W,r}$.
Let $J: \bar \UU\ra B_r(f)$, $J(g)= g-f|\, W$ denote this embedding.
Since $\HH_0(\eps/2)$ is compact, $J(\bar \UU)$ is compact as well.

On the other hand, by Property (P) and continuity of $P$,
there is a $\rho>0$ such that $P(j_{f,W,\rho}(B_\rho))\subset \UU$.
Hence $P_f: = J\circ P \circ j_{f,W,\rho}: B_\rho\ra B_r$ is a compact
Banach holomorphic retraction.

Let us consider its real-symmetric part
$P_f^\R: B^\R_\rho \to B_r^\R $.
It is a compact real analytic Banach retraction,
so by  Lemma \ref {dp}, 
the set  $\Fix P^\R_f$  near $f$ is a real analytic finite-dimensional
submanifold of $B_\rho$.
But since $\ZZ^\R$ is compact, topology induced on it from the Banach ball
coincides with its own topology (induced from the whole space $\HH_0$).
Hence $\ZZ^\R$ is a finite-dimensional topological manifold near $f$.
Since $f\in \ZZ^\R$ is arbitrary,
the first step of the proof is completed.

\msk
By Lemma \ref {contractible}, the space  $\EE^\R$ is contractible.
Since $\HH_0^\R$ is homeomorphic to it, it is also contractible.
Since $\ZZ^\R$ is a retract of $\HH_0^\R$, it is contractible  as well.
(If $h_t$ is a homotopy that contracts $\HH_0^\R$ to a point,
then $P\circ h_t : \ZZ^\R \ra \ZZ^\R$ does the same to $\ZZ^\R$.)
But  the only contractible compact finite-dimensional  manifold (without
boundary) is a point
(since otherwise the top homology group $H^n(M)$ is non-trivial.)
This concludes the 2nd step.

\msk
Thus,  $P_f(B_\rho^\R)=\{0\}$.
Since $P_f: B_\rho\ra B_r$ is holomorphic (as a Banach map), $ P_f
(B_\rho)=\{0\}$.

Let us show that a small neighborhood $\VV$  of $f$ in $\ZZ$
belongs  to the neighborhood $\UU\subset \HH_0(\eps/2)$ considered above.
Otherwise  $f\in \cl(\ZZ\sm \UU)$.  Since the notion of closedness in $\HH_0$ is
given in terms of sequences,
there would exist a sequence $f_n\in \ZZ\sm \UU$ converging to $f$.
By Property (P), the maps $f_n=Pf_n$ would eventually belong to
$\HH_0(\eps/2)$,
and hence to the neighborhood  $\UU$ -- contradiction.

Thus,  we have   $\VV\subset \UU\subset B_r$.
Shrinking $\VV$ if needed, we make $J(\VV)\subset B_\rho$  and hence
$$
  J(\VV)=J(P(\VV)) \subset P_f(B_\rho)=\{0\}.
$$
Since $J$ is injective, $\VV=\{f\}$.
Thus,  $f$ is an isolated point in $\ZZ$.
But since $\ZZ= P(\HH_0)$ is connected, we conclude that $\ZZ=\{f\}$,
which completes the last step of the proof.
\qed

\comm{
Let $\psi_{f,n}:\inter K(f) \to \C$

$(\lambda,f_\lambda(z))$ provides a function $\D \times
\D_{1/64} \to \D \times \D_{1/64}$ for every $n \geq 1$, which is
holomorphicBy iteration, we see that $(\lambda

Thus for every $n \geq 1$,
$f_\lambda^n$ restricts to a holomorphic function $\D_{1/64} \to \D_{1/64}$.

and since it
depends holomorphicall
It follows that $(\lambda,z) \mapsto f_\lambda^n(z)

Since $I \circ \gamma$ is holomorphic,
$(\lambda,z) \mapsto \gamma_\lambda^n(z)$, and hence
$(\lambda,z) \mapsto \psi_{\gamma_\lambda,n}(z)$, are holomorphic maps from
$\D \times \D_{1/16}$ to the unit disk,
for every $n \geq 1$.

By (\ref {}), $\D_{1/4} \subset K(f) \subset W$, so
$w_\lambda(z)$, and hence $f_\lambda(z)=w_\lambda(z)+f(z)$,
is a holomorphic function of $(\lambda,z) \in \D \times \D_{1/4}$.

Fix $0<R<1/4$.  Letting $j:B_{W,r_0} \to \BB_{$ be the restriction of
$w \mapsto w(z)+f(z)-z^d$. If $\gamma$ is holomorphic then $j \circ \gamma=
I \circ j_{f,W,r_0} \circ \gamma$ is holomorphic, so by Lemma
\ref {crit}
$j_{f,W,r_0} \circ \gamma$ (which is continuous since $\gamma$ is
holomorphic and $j_{f,W,r_0}$ is continuous) is path holomorphic.

\end{pf}
}
}

\section{Almost periodicity and retractions}

We now turn to the dynamical construction of retractions.
The presence of enough compactness, together with the non-expansion of the
Carath\'eodory metric, allows us to implement the
notion of {\it Almost Periodicity} (adapted appropriately to the
cocycle setting).

\subsection{Almost periodic cocycles}\label{cocycles}

We will now discuss cocycles with values in a Hausdorff
topological semi-group $\SSS$.
Since we aim to eventually take
$\SSS=\Hol^\R(\HH_0,\HH_0)$, we allow for the possibility that $\SSS$ is
not metrizable,  neither satisfy the First Countability Axiom, however we
will always assume that $\SSS$ is {\it sequential}
in the sense that
the notions of continuity, closedness and compactness
can be defined in terms of sequences.

\def\SSS{{\mathcal{S}}}


\comm{
Let $Q=\{(m,n)\in \N\times \N:\ n>m\}$.            
An $\SSS$-{\it cocycle} is a map $g:  Q \ra \SSS$, $(m,n)\mapsto g^{mn}$,  such that
$g^{mn} g^{lm}  = g^{ln}$. 
Letting $f_n:= g_{n, n+1} \in \SSS$, we obtain
\begin{equation} 
    g^{mn}= f_{n-1} \dots f_{m+1}  f_m,
\end{equation}
and vice versa, any sequence $f_n\in \SSS$ determines a cocycle by means of (\ref{cocycle}). 
}

A {\it subcocycle} is the restriction $G^{k_m,k_n}$ of a cocycle $G$ to a subsequence of $\N$.
More formally, it is a pullback of $G$ under a strictly monotone embedding $k: \N\ra \N$. 



The $\om$-limit set of a cocycle, $\om(G)$, is the set of all existing 
 $\lim G^{m,n}$ as $m\to \infty$ and $n-m\to \infty$. 

A cocycle is called {\it almost periodic}
if the family $\{G^{m,n} \}_{(m,n)\in Q}$ is precompact in $\SSS$.
The $\om$-limit set of an almost periodic cocycle is compact.

We say that a cocycle {\it converges} if there is the limit 
$
   G^{m,\infty} : = \lim_{n\to \infty} G^{m,n}
$
for every $m \in \N$.  The cocyclic rule extends to the limits of converging cocycles:
$$
    G^{m,\infty} = G^{n,\infty}\,  G^{m,n}, \quad (m,n)\in Q.
$$ 

We say that a cocycle {\it double converges} if there exists the limit $G^{\infty\infty}:= \lim_{m\to \infty}
G^{m,\infty}$.
The cocyclic rule extends to the limits of double converging cocycles: 
$$
           G^{m,\infty} = G^{\infty,\infty}\,  G^{m,\infty},\quad
G^{\infty,\infty} =
(G^{\infty,\infty})^2.
$$
In particular, $G^{\infty,\infty}$ is an {\it idempotent}.

\begin{lem}
  An almost periodic cocycle has a double converging subcocycle.
\end{lem}

\begin{pf}
  A converging subcocycle is extracted by means of the diagonal process.
Selecting then a converging subsequence of the $G^{m,\infty}$, we obtain a double converging subcocycle. 
\end{pf}

\comm{****
\begin{lem}
  The double limit $P:=g^{\infty\infty}$ of a double converging cocycle is an idempotent.
\end{lem} 

\begin{pf}
  For any $m<n$,  we have a cocyclic rule: $g_{m\infty} = g_{mn} g_{n\infty}$. 
Letting $n\to \infty$, we obtain $g_{m\infty} = g_{m\infty} P$.
Letting now $m\to\infty$, we obtain $P^2=P$.    
\end{pf}
******}

\begin{cor}\label{idempotent in om}
  The $\om$-limit set of an almost periodic cocycle contains an idempotent. 
\end{cor}

We endow the space of cocycles with the pointwise convergence topology:
$$
G_k\to G\quad \mathrm{if}\quad G^{m,n}_k \to G^{m,n}\quad \mathrm{for\ all} \ (m,n)\in Q.  
$$

The {\it shift} $T$ in the space of cocycles is induced  by the embedding $\N\to \N$, 
$n\mapsto n+1$. In other words, $(T G)^{m,n} = G^{m+1, n+1}$, $(m,n)\in Q$.

If $G$ is almost periodic then all its translates $\{T^n G\}_{n=0}^\infty$ form a precompact family of cocycles.  

Given a function $\rho: \SSS\ra \R_{\geq 0}$, 
we say that a cocycle is {\it uniformly $\rho$-contracting} 
 if for any $\gamma>0$ there exists an $N$
such that $\rho (G^{m,n})< \gamma $ for any $(m,n) \in Q$ with $n-m\geq N$. 

A continuous function $\rho: \SSS\ra \R_{\geq 0}$ is called {\it Lyapunov} 
$$
    \rho (F_l F F_r) \leq \rho(F)\quad {\mathrm {for\ any}}\ F,F_l,F_r\in \SSS. 
$$ 
The next assertion will not be directly used  but can serve as a model for what follows:

\begin{prop}
  Let $G$ be an almost periodic cocycle that has a Lyapunov function $\rho$. 
If $\rho(e)=0$ for any limit idempotent $e\in \om(G)$ 
then the cocycle is uniformly $\rho$-contracting. 
\end{prop}

We leave it as an exercise.

\comm{****
\begin{pf}
  Otherwise there exists a $\gamma>0$ and two non-decreasing  sequences $q_k\in \N$ and $n_k\to \infty$ such that
$\rho( g_{m_k}^{q_k+n_k}) \geq \gamma$.  Since $g$ is almost periodic, the sequence of cocycles $T^{q_k} g$ admits a
converging subsequence. Let $h$ be a limit cocycle. Then
$
   \rho( h_{0n}) \geq \gamma>0 \quad \mathrm{for \ any}\ n>0.
$ 
Since $\rho$ is Lyapunov, we have $\rho( h_{lm}) \geq \gamma$ for all $(l,m) \in Q$. 
Hence $\rho( \phi) \geq \gamma$ for all $\phi \in \om (h)$. 
In particular, $\rho(e)>0$ for any idempotent $e\in \om(h)$ from  Lemma \ref{idempotent in om}.
Since $\om(h)\subset \om (g)$, we arrive at a contradiction. 
\end{pf}
************}

We will need a more general form of the above proposition.
Assume that we have two continuous functions,  $\rho'\geq \rho\geq 0 $ on $ \SSS $,
which are not assumed to be individually Lyapunov, 
but rather possess a joint Lyapunov property (adapted to the cocycle):
\begin{equation}\label{joint Lyap}
   \rho(G^{l,n}) \leq \rho(G^{l,m}) \quad  
 \mathrm{and}\quad   \rho(G^{l,n}) \leq \rho'(G^{m,n}) \quad {\mathrm {for\ any}}\  l < m  < n. 
\end{equation}
We call it a {\it Lyapunov pair} for the cocycle. 

\comm{****
\begin{lem}\label{contraction with two norms}
 Let $g$ be an almost periodic cocycle that has a Lyapunov pair $(\rho, \rho')$. 
If $\rho'(e)=0$ for any limit idempotent $e\in \om(g)$  then $g$ is $\rho$-contracting. 
\end{lem}

\begin{pf}
  Otherwise there exists a $\gamma>0$ and a non-decreasing  sequence $m_k\in \N$ such that
$\rho( g_{m_k}^{m_k+ k}) \geq \gamma$.  Since $g$ is almost periodic, the sequence of
cocycles $T^{m_k} g$ admits a
converging subsequence. Let $h$ be a limit cocycle. Then
$\rho( h_0^n) \geq \gamma>0$ for any $n>0$ (by continuity of $\rho$ ane the first part of (\ref{joint Lyap})).   
By the joint Lyapunov property, we have $\rho'( h_l^m) \geq \gamma$ for all $(l,m) \in Q$. 
Hence $\rho' ( \phi) \geq \gamma$ for all $\phi \in \om (h)$. 
In particular, $\rho'(e)>0$ for any idempotent $e\in \om(h)$ from  Lemma \ref{idempotent in om}.
Since $\om(h)\subset \om (g)$, we arrive at a contradiction. 
\end{pf}
*****}

We will also need a uniform version of the above lemma, over a family of cocycles.
Let $\GG$ be a family of cocycles $G_s$ labeled by an element $s$ of some set
$\Sigma$.  We say that $\GG$ is {\it uniformly almost periodic} if the whole family
$G^{m,n}_s$,  $s\in \Sigma$, $(m,n)\in Q$,
is precompact in $\SSS$.  Then $\om(\GG)\subset \SSS$ stands for the set of the limits of all converging sequences
$G^{m_k, n_k}_{s_k}$ as $m_k\to \infty$ and  $n_k-m_k\to \infty$. 

We say that the family is {\it uniformly $\rho$-contracting} if 
for any $\gamma>0$ there exists an $N$
such that $\rho (G^{m,n}_s)< \gamma $ for any $s\in \Sigma$, $m\in \N$ and $n\geq m+ N$.

\begin{lem}\label{contraction with two norms}
Let $\GG$ be a uniformly almost periodic family of cocycles.
Let $(\rho, \rho')$ be  a Lyapunov pair  for all cocycles in $\GG$. 
If $\rho'(e)=0$ for any limit idempotent $e\in \om(\GG)$  then $\GG$ is uniformly $\rho$-contracting. 

\end{lem}

\begin{pf}
  Otherwise there exists a $\gamma>0$, a sequence $s_k\in \Sigma$
and two non-decreasing  sequences $q_k\in \N$ and $n_k\to \infty$ such that
$\rho( G^{q_k, q_k+ n_k}_{s_k}) \geq \gamma$.  Since $\GG$ is uniformly almost periodic, the sequence of
cocycles $T^{q_k}G_{s_k}$ admits a
converging subsequence. Let $G$ be a limit cocycle. Then
$\rho( G^{0,n}) \geq \gamma>0$ for any $n >0$ (by continuity of $\rho$ and the first part of (\ref{joint Lyap})).   
By the second part of (\ref{joint Lyap}), we have $\rho'( G^{m,n}) \geq \gamma$ for all $(m,n) \in Q$. 
Hence $\rho' ( \phi) \geq \gamma$ for all $\phi \in \om (G)$. 
In particular, $\rho'(G)>0$ for any idempotent $e\in \om(G)$ from  Lemma \ref{idempotent in om}.
Since $\om(G)\subset \om (\GG)$, we arrive at a contradiction. 
\end{pf}


\subsection{Tame spaces}\label{tame spaces}
%

Let $X$ be a (sequential)  topological space endowed with a
continuous metric $d: X \times X\ra \R_{\geq 0}$
that is compatible with the
topology on compact subsets of $X$ (but not necessarily  on $X$).
We say that $X$ is {\it tame}
if the following properties hold:
\begin{enumerate}
\item There exists a filtration of compact subsets,  $X_1\subset X_2\subset\dots \subset X$ such that $\cup X_i=X$;
\item Each compact set in $X$ is contained in some $X_i$;
\item a set is open in $X$
if and only if its intersection with any compact subset of $X$ is relatively open.
\end{enumerate}

A family of continuous
maps $F_s:X \to X'$, $s\in \Sigma$,  between tame spaces (with metrics $d$ and $d'$ respectively)
is called {\it equicompact}
if for every compact set $K \subset X$ there exists a compact set
$K' \subset X'$ such that $F_s(K) \subset K'$ for all $s\in \Sigma$.

\ssk \nin $\bullet$ 
An equicompact family $\{F_s\}$ is called {\it equicontinuous on compact sets}
 if  for every compact set $K \subset X$, 
for every $x \in K$ and $\epsilon>0$ there exists
$\delta=\delta(K,x,\epsilon)>0$\footnote {In fact, by compactness of $K$,
one may take $\delta=\delta(K,\eps)$ here.}
such that for every $s \in \Sigma$, we have
$$
\mbox{  if $d(x,y)<\delta$ and $y \in K$ then $d'(F_s(x),F_s(y))<\epsilon$.}
$$

\ssk \nin $\bullet$
A sequence $\{F_n\}$ is called {\it  uniformly converging on compact sets}
if $\{F_n\}$ is equicompact and there exists a continuous map $F:X \to X'$
such that for every compact set $K \subset X$, we have
$$ 
  \lim_{n \to \infty} \sup_{x \in K} d' (F_n(x),F(x))=0.
$$  
Notice that in 
this case $F_n$ is necessarily equicontinuous on compact sets.

\ssk \nin $\bullet$
A family of cocycles $G_s$, $s\in \Sigma$, is called 
 {\it uniformly contracting on compact sets}
if for any compact subset $K\subset X$ and any $\gamma>0$,
there exists a compact set $K'$ and an $N$ such that
$$
G^{m,n}_s(K)\subset K'\quad \mathrm{and}\quad \mbox{$\diam (G^{m,n}_s (K)) < \gamma$ for all
$s\in \Sigma$, $m\in\N$ and  $n\geq m+N$.}   
$$

The space of continuous maps $F: X\ra X'$ between two tame spaces is endowed with topology of
uniform convergence on compact subsets.

\begin{lemma} \label {comp}

A sequence of maps $F_n :X \to X'$ between tame spaces is precompact
if and only if it is equicontinuous on compact sets. 

\end{lemma}

\begin{pf}
In the non-trivial direction, it follows from 
the Ascoli-Arzela's Theorem on each $X_i$ and the diagonal argument. 
\end{pf}

Given a tame space $X$, 
let $\SSS\equiv \SSS_X$ be the topological semigroup of all continuous
weak contractions of $X$
(endowed with topology of uniform convergence on compact subsets). 
Idempotents in this semigroup are retractions.   

\begin{lem}\label{uniform cotraction in tame spaces}
 Let $X$ be a tame space, and    
let $\GG$ be  a uniformly almost periodic family of cocycles
$G^{m,n}_s$, $s\in \Sigma$,
with values in the semigroup $\SSS_X$. 
If all limit retractions $P\in \om (\GG)$ are constants, then $\GG$ is uniformly contracting
on compact sets. 
\end{lem}

\begin{pf}
   Since $\GG$ is uniformly almost periodic, the family of maps $G^{m,n}_s$  is  equicompact: 
for any $i\in \N$ there exists $j=j(i)\geq i$ such that 
$G^{m,n}_s (X_i)\subset X_j$ for all $s\in \Sigma$,  $(m,n)\in Q$. 

Let $\rho_i: \SSS\ra \R_{\geq 0}$ be defined as $\rho_i(G)= \diam G(X_i)$.
Obviously, these functions are continuous and form a monotonically increasing sequence. 
Moreover, for any $i\in \N$, the functions $\rho:=\rho_i$ and  $\rho':=\rho_{j(i)}$
form a Lyapunov pair for any cocycle $G \in \GG$.  Indeed, 
for any  $l< m< n$ we have:
$$
    \rho (G^{l,n}) = \diam G^{l,n}(X_i) \leq \diam G^{l,m} (X_i)=\rho(G^{l,m})
$$
(where the estimate holds since the semigroup $\SSS$ consists of weakly contracting maps)
$$
   \rho (G^{l,n}) =  \diam G^{l,n}(X_i) \leq  \diam G^{m,n} (X_j)  =
\rho'(G^{m,n}). 
$$
(where the estimate holds since $G^{l,m}(X_i)\subset X_j$).  

Since all retractions  $P\in \om  (\GG)$ are constants, we have $\rho'(P)=0$ for any of them. 
By Lemma~\ref{contraction with two norms},
the family $\GG$ of cocycles is uniformly $\rho$-contracting, so 
for any $\gamma>0$, there exists an $N$ such that 
$\diam G^{m,n}_s(X_i) <\gamma$ as long as $n\geq m+N$.
  
 Since $i$ is arbitrary, we are done.
\end{pf}

\subsection{Proof of Theorem \ref {trivial}}

Notice that $\HH_0$, with the Carath\'eodory metric,
is tame in the sense of \S \ref {tame spaces} 
(take, e.g., $X_i: = \HH_0(2^{-i})$ as a filtration).
By the Schwarz Lemma, 
$\Hol(\HH_0,\HH_0)$ is a sub-semigroup of $\SSS_{\HH_0}$, which turns out to
be closed:

\begin{lemma} \label {pa}

If $F_n \in \Hol(\HH_0,\HH_0)$ converges uniformly on compact sets to  a map
$F$ then $F \in \Hol(\HH_0,\HH_0)$.

\end{lemma}

\begin{pf}

We have to show that if $\gamma:\D \to \HH_0$ is a holomorphic path
then $F \circ \gamma$ is a holomorphic path as well.
Let $0<\rho<1/4$, and let $I_\rho: \HH_0\ra \BB_{\D_\rho}$ be the
restriction operator (c.f. (\ref {1/4})).
The sequence of maps $\{I_\rho \circ F_n \circ \gamma\}_n$ converges uniformly on compact sets
to $I_\rho \circ F \circ \gamma$. By Lemma \ref {crit}, each $I_\rho \circ F_n \circ \gamma$ is
holomorphic in the usual Banach sense, so the limit
$I_\rho \circ F \circ \gamma$ is holomorphic as well.  
By Lemma \ref {crit}, $F \circ \gamma$ is path holomorphic.
\end{pf}

Property (H1) implies that the family $\GG$ is uniformly almost periodic.
If (C1) does not hold then by Lemma \ref {uniform cotraction in tame spaces},
there exists a
sequence $G_k \in \GG$ and $(m_k,n_k)\in Q$ with $n_k-m_k \to \infty$,
such that $G^{m_k,n_k}_k$
converges uniformly on compact sets to a non-constant retraction $P \in
\SSS_{\HH_0}$.  By Lemma \ref {pa}, $P$ is path holomorphic, concluding the
proof.
\qed

\section{Horseshoe}

\subsection{Beau bounds and rigidity yield the horseshoe} 

\subsubsection{Complex case}
Let $\FF$ be a family of disjoint little Multibrot sets $\MM_k$
(encoding certain renormalization combinatorics).
We say that {\it beau bounds are valid for} $\FF$, if they are valid for the
family of infinitely renormalizable maps $f$ 
whose renormalizations $R^n f$ have combinatorics $\bar \MM = (\MM_n)_{n=0}^\infty$ with
$\MM_n\equiv \MM_{k_n}\in \FF$. (We will loosely say that ``$\bar \MM$ is in $\FF$''.) 

We say that $\FF$ {\it is rigid} if for any combinatorics $\bar \MM$ in  $\FF$,  
there exists a unique polynomial $p_c: z\mapsto z^d+c$ which is infinitely renormalizable with this combinatorics. 

\begin{rem}
It is conjectured that any family $\FF$ is in fact rigid
(which would imply that the Multibrot set is locally connected at all infinitely renormalizable parameter values, 
and hence would prove MLC, for all unicritical families).  
\end{rem}

In \cite {Davoud}, it is shown that beau bounds implies rigidity for a large class of combinatorics, 
and in fact for all combinatorics for which beau bounds have been proved (\cite {K}, \cite {KL1}, \cite {KL2}).

\begin{rem}
  In the quadratic case, this rigidity result had been established in \cite{puzzle}.
The  work \cite {Davoud} makes use of recent advances:
a new version of the {\it Pullback argument} developed in \cite{AKLS}.
\end{rem}

A {\it semi-conjugacy} between two dynamical systems $F: X\ra X$ and $G: Y\ra Y$ is a continuous {\it surjection}
$h: X\ra Y$ such that $h\circ F = G\circ h$. 
  
Let  $\Sigma\equiv \Sigma_\FF= \FF^\Z$ be the symbolic space with symbols from $\FF$, 
and let  $\sigma: \Sigma\ra \Sigma$ be the corresponding two-sided shift.
A map $h:\Sigma \to \CC$ is called {\it combinatorially faithful}
if for any combinatorics  $\bar \MM = (\MM_n)_{n=-\infty}^\infty \in \Sigma$, 
the image $f=h(\bar \MM)$ is renormalizable with combinatorics $\MM_0$,  
and $h$ semi-conjugates $\sigma$ and $R|\AA$.

\begin{thm}\label{complex horseshoe}
Assume a family $\FF$ has beau bounds and is rigid.
Then  there exists a precompact $R$-invariant set $\AA\subset \CC$ 
and a combinatorially faithful semi-conjugacy $h: \Sigma \ra \AA$. 
Moreover, $R$ is exponentially contracting along the leaves of the hybrid
lamination of $\AA$, with respect to the
Carath\'eodory metric.\footnote {It follows from rigidity that the hybrid
lamination of $\AA$ consists of all infinitely renormalizable maps whose
renormalizations have combinatorics in $\FF$.}
\end{thm}

\begin{pf}
Since $\FF$ has beau bounds,
there exists $\eps>0$ such that 
$$
 \mod R^n p \geq \eps, \quad n=0,1,\dots
$$  
for any polynomial $p: z\mapsto z^d+c$
which is infinitely renormalizable with combinatorics $\FF$.  

Let $\overline\MM = (\MM_n)_{n \in \Z} \in \Sigma$. 
For any $-n\in \Z_-$, 
there exists a polynomial $p_{c_n}$ 
which is infinitely renormalizable with combinatorics $(\MM_{-n}, \MM_{-n+1},\dots)$. 
Then for any $l\in \Z$, $l\geq -n$,
the germ $f_{n,l}:= R^{n+l} p_{c_{-n}}\in \CC(\eps)$ is infinitely renormalizable with combinatorics $(\MM_l, \MM_{l+1},\dots)$.
As this family of germs is precompact,   for any $l$ we can select a subsequence $f_{n(i),l}$
converging to some $f_l\in \CC(\eps)$ as $n(i)\to -\infty$. 
This map is infinitely renormalizable with combinatorics  $(\MM_l,\MM_{l+1},\dots)$.%
\footnote{This actually needs a little check-up: compare the argument four paragraphs down.}
Using the diagonal procedure (going backwards in $l$) we ensure that $R f_{l-1} = f_l$. 
Thus, we obtain a bi-infinite sequence of maps $f_l\in \CC(\eps)$ such that  
$f_l$ is renormalizable with combinatorics $\MM_l$ and  $R f_{l-1} = f_l$.

Assume there exist two such sequences, $(f_l)_{l\in \Z}$ and $(\tl f_l)_{l\in \Z}$.
Since $\FF$ is rigid, the hybrid class of any $f_l$ is uniquely determined by the 
renormalization combinatorics $(\MM_l, \MM_{l+1}, \dots )$, so for any $l\in \Z$, the germs $f_l$ and $\tl f_l$ are hybrid equivalent.  
But by Theorem \ref{contr}, 
the renormalization is exponentially contracting with respect to the Carath\'eodory metric in the hybrid lamination. 
Hence there exist $C>0$ and $\la\in (0,1)$ such that
\begin{equation}\label{exp close}
   d_{\chi(f_0)} (f_0, \tl f_0) \leq C(\epsilon) \la^n d_{\chi(f_n)} (f_{-n}, \tl f_{-n}).
\end{equation}
Letting $n\to \infty$ we see that $f_0=\tl f_0$.
For the same reason, $f_l=\tl f_l$ for any $l\in \Z$.

Thus, we obtain a well defined equivariant map $h: \overline\MM \mapsto f$, 
where  $f\equiv f_0$ is a polynomial-like germ for which  there exists a bi-infinite
sequence $f_l\in \CC(\eps)$, $l\in \Z$, such that $f_l$ is renormalizable with combinatorics $\MM_l$. 
Let $\AA\subset \CC$ consist of all such germs, which makes $h$ surjective by definition.    

To see that $h$ is continuous, consider a sequence $\overline \MM^{(k)} \to
\overline \MM$ in $\Sigma$, and let $f^{(k)}_l=h(\sigma^l(\overline
\MM^{(k)}))$.  We need to show that $h(\overline \MM^{(k)}) \to h(\overline
\MM_k)$.  By passing through an arbitrary
subsequence, we may assume that for each $l
\in \Z$, $f^{(k)}_l$ converges in $\CC(\epsilon)$ to some $f_l$.  By
definition of convergence in $\Sigma$, for each $l \in \Z$ and for each $k$
sufficiently large, $f^{(k)}_l$ is renormalizable with combinatorics in
$\MM_l$, i.e., $\chi(f^{(k)}_l) \in \MM_l$.

Let us show that $f_l$ is renormalizable with combinatorics $\MM_l$.
If $\MM_l$ is a primitive copy, then it is closed, which readily implies
that $\chi(f_l) \in \MM_l$.  If $\MM_l$ is a satellite copy, its closure is
obtained by adding the root, so we also need to guarantee that
$\chi(f^{(k)}_l)$ does not converge to the root.  But for $k$ large,
$\chi(f^{(k)}_l)$ belongs to a subcopy of $\MM_l$ (consisting of those
polynomials in $\MM_l$ whose renormalization has combinatorics $\MM_{l+1}$),
which is at definite distance from the root of $\MM_l$, so we can again
conclude that $\chi(f_l) \in \MM_l$.

By continuity of any renormalization operator with fixed combinatorics,
we also conclude $R f_l=f_{l+1}$.  By the definition of $h$, $h(\overline
\MM)=f_0=\lim f^{(k)}_0=\lim h(\overline\MM^{(k)})$, as desired.
\end{pf}

\subsubsection{Real horseshoe}
Let  $\FF^\R$ stand for the family of all real renormalization combinatorics with minimal periods.
Let $\si : \Sigma \to \Sigma$ be the corresponding shift.
While neither beau bounds, nor (complex) rigidity have been established for $\FF^\R$,%
\footnote {They have been established, however, 
for the family of primitive real combinatorics \cite{KL1}, 
which covers all real combinatorics except period doubling.}
the proven beau bounds and rigidity for  {\it real-symmetric germs} is enough to construct the renormalization horseshoe. 
Moreover, we will show that in this case, the combinatorially faithful semi-conjugacy $h$ is actually a homeomorphism:

\begin{thm} \label {real horseshoe}
  For the family $\FF^\R$, 
  there is an $\R$-invariant set $\AA\subset \CC^\R$ 
and a combinatorially faithful homeomorphism
$h: \Sigma \ra \AA$. 
Moreover, $R$ is exponentially contracting along the leaves of
the hybrid lamination of $\AA$ (endowed with the Carath\'eodory metric),
which contains all infinitely renormalizable real-symmetric germs.  
\end{thm}

\begin{pf}

The construction of the horseshoe $\AA$, along with the combinatorially faithful semi-conjugacy  $h$,
is basically the same as in Theorem \ref{complex horseshoe}, with the following adjustments:
\begin{enumerate}
\item  The polynomials $p_{c_n}$ should be selected to be real, $c_n \in \R$;
\item Theorem \ref {real beau} ({\it beau bounds for real maps}) 
       provides the needed compactness for the construction of the maps $f_l$;
\item Rigidity of $\FF$ is replaced with {\it rigidity for real polynomials}, obtained in \cite {puzzle}, \cite {GS}
(quadratic case) and in \cite {KSS} (arbitrary degree):%
\footnote{This result also follows from combination of 
\cite{AKLS} (dealing with at most finitely renormalizable maps)
and \cite{Davoud} (dealing with the infinitely renormalizable situation).}
any real renormalization combinatorics determines a single {\it real-symmetric} hybrid leaf;
\item {\it  Exponential contraction} is obtained by combining Theorems \ref {beau
complex} and \ref {contr}.
\end{enumerate}
The injectivity of $h$ follows from the injectivity of the renormalization
operator acting on real p-l maps \cite[p. 440]{MvS}.
What is left, is  to verify continuity of $h^{-1}$.  Since convergence in
$\Sigma$ means coordinatewise convergence, it is equivalent to
the following statement:

\begin{lem}\label{bounded comb for each n}
Let  $(\bar \MM^{(j)})_{j=1}^\infty$ be a sequence of symbolic strings 
$\bar \MM^{(j)}=(\MM^{(j)}_n)_{n\in \Z}$ in $\Sigma$ 
such that the corresponding  germs $f_j\equiv h(\bar \MM^{(j)})$ converge to some
$f_\infty\in \AA $.
Then for any $n\in \Z$, the combinatorics $\MM^{(j)}_n$
of $R^n f_j$ eventually coincides with that of $R^n(f_\infty)$.
\end{lem}

Notice that Lemma \ref {bounded comb for each n} is clear for $n=0$: 
a real perturbation of  a twice renormalizable real map 
is renormalizable (at least once) with the same combinatorics 
(since on the boundary of the renormalization windows the maps are not twice renormalizable, 
see \S \ref{renorm sec}).  Since the renormalization operator acts continuously on
$\AA$, we have  for any $n\geq 0$ that $R^n f_j\to R^nf_\infty$ as $j\to \infty$. 
It follows that Lemma \ref {bounded comb for each n} holds for all $n\geq 0$ as well.

In order to prove it inductively for $n<0$, 
it is enough to show that the hypothesis of Lemma \ref {bounded comb for each n} imply
that $R^{-1} f_j \to R^{-1} f_\infty$.  
To this end, it is sufficient to prove that 
the {\it  renormalization combinatorics of the germs $R^{-1} f_j$ are bounded}.  
Indeed, in this case any limit  $g$  of these germs is renormalizable and $Rg=f_\infty$. 
By injectivity of the renormalization operator, $g= R^{-1} f_\infty$, 
and the conclusion follows.

Boundedness of the renormalization combinatorics 
follows from an analysis of the domain of analyticity  of
limits of renormalized germs:

\begin{lem}\label{Dichotomy}
Let $\tilde f_j \in \CC^\R$, $j \geq 1$,
be a sequence of renormalizable germs. 
If the renormalization periods of the $\tilde f_j$ are going to infinity,     
then any limit of renormalizations $R \tilde f_j$  is either 
a unicritical polynomial or its real trace
has a bounded domain of analyticity.
\end{lem}

See Appendix \ref {analytic} for a proof (unlike the previous parts of this paper,
it relies on the fine combinatorial and
geometric structure of {\it one} renormalization).

Let us apply Lemma \ref {Dichotomy} to $\tilde f_j=R^{-1} f_j$.
Notice that $f_\infty=\lim R \tilde f_j$ cannot be a unicritical polynomial
since those are never anti-renormalizable.
It can neither have a bounded domain of analyticity, due to McMullen's result \cite{towers} implying that
the real trace of
$f_\infty$ (which is {\it infinitely} anti-renormalizable with {\it a priori} bounds)
extends analytically to $\R$.
So, both options offered by  Lemma \ref {Dichotomy} are impossible in our situation,
and hence the renormalization periods of the germs $\tl f_j$ must be bounded. 
This concludes the proof of Lemma \ref {bounded comb for each n}, and thus of
Theorem \ref {real horseshoe}.
\end{pf}



\appendix

\section{Analytic continuation of the first renormalization} \label {analytic}

\subsection{Principal nest and scaling factors}

Through this section, we consider a renormalizable unimodal map $f:I \to I$, of period $p$,
with a polynomial-like extension in $\CC^\R(\epsilon_0)$
for some fixed $\epsilon_0>0$.  For simplicity of notation, we will also assume that $f$
is even.  We will also assume that we can write $f(x)=\psi(x^d)$ for some
diffeomorphism $\psi$ with non-positive Schwarzian derivative.  All arguments
below can be carried out without the extra assumptions with only technical
changes, but in the situation arising in our application ($f$ is infinitely
anti-renormalizable with {\it a priori} bounds) they are indeed
automatically satisfied.

Below $C>1$ stands for a constant which may only depend on $\epsilon_0$.

Recall that a closed interval $T \subset I$ which is symmetric
(i.e., $f(\partial T)$ is a single point) is called {\it nice} if
$f^k(\partial T) \cap \inter T=\emptyset$, $k \geq 1$.  If 
the critical point returns to the interior of a nice interval $T$, then
we let $T'$ be the central component of the first return map to $T$.  We let
$\lambda(T)=|T'|/|T|$ be the {\it scaling factor}.

We define the {\it principal nest} $I_n$, $n \geq 0$ as follows.  Since $f$ is
renormalizable, it has a unique orientation reversing fixed point $p$.  
Its preimage $\{p,-p\}$ bounds a nice interval, which we denote $I_0$. 
Then we define $I_{n+1}$ inductively as $I_n'$, i.e., $I_{n+1}$ is
the central component of the first return map to $I_n$.

Let us assume from now on that $p>2$. 
Under this condition, $I_{n+1} \ssubset \inter I_n$ for every $n \geq 0$.  
Let $\lambda_n=|I_{n+1}|/|I_n|$ be the corresponding scaling factors.

Let $g_n$ be the first return map to $I_n$.  
We say that $g_n$ is {\it central} if $g_n(0) \in I_{n+1}$.  
We define a sequence $(j_k)_{k \geq 0}$ inductively so that
$j_0=0$ and $j_{k+1}$ is the minimum $n>j_k$ such that $g_{n-1}$ is non-central.  
Since $f$ is renormalizable, $g_n$ is central  for all sufficiently large $n$, 
so the sequence $(j_k)$ terminates at some $N=j_{\kappa}$.
We call $\kappa$ the {\it height} of $f$.

We have the following basic estimates on the scaling factors (see \cite{Ma1}): 
 
 
\proclaim {A priori bounds.}
{
We have:  $\lambda_{j_k} \leq 1-C^{-1}$ for every $k$.
Moreover. the maps $g_{j_k}: I_{j_{k+1}}\ra I_{j_k}$ are compositions of power maps $x\mapsto x^d$ 
and diffeomorphisms with bounded distortion. 
} 

\begin{cor}
\item $\lambda_{n+1} \leq C \lambda_n^{1/d}$,
\item $\lambda_{j_{k+1}} \leq C \prod_{n=j_k}^{j_{k+1}-1} \lambda_n^{1/d}$.
\end{cor}

Let $v_n$ be the {\it principal return times}, i.e. $f^{v_n}|I_{n+1}=g_n$.  
In particular, $v_N=p$ is the renormalization period.  
We let $g=f^p$.  
Then $g: J \to J$ is the {\it unimodal pre-renormalization} of $f$,
where $ J=\cap_{n \geq 0} I_n$.


\begin{lemma} \label {I'_0}

We have $T_p \subset I_{\max \{0,N-1\}}$.

\end{lemma}

\begin{pf}

If $\kappa=0$ then clearly $T_p=I_1 \subset I_0$.

Assume that $\kappa \geq 1$ and hence $N \geq 1$.
The interval $T_p$ is the smallest interval containing $0$ whose boundary
is taken by $f^p$ to the boundary of $I_0$.  
Hence it is enough to show that $f^p (\partial I_{N-1}) \subset \partial I_0$.

It is easy to see that for $k \geq 1$ we have  
$$
f^{s v_{j_{k-1}}}(0) \in I_{j_k-s} \setminus I_{j_k+1-s}\quad \text{for}\quad 1 \leq s \leq j_k-j_{k-1}.
$$  
We conclude that      
$$
   v_{j_k} \geq v_{j_{k-1}-1} + (j_k-j_{k-1}) v_{j_{k-1}} = \sum_{n=j_{k-1}-1}^{j_k-1} v_n
=v_{j_{k-1}}+ \sum_{n=j_{k-1}-1}^{j_k-2} v_n, \quad \text{for}\quad j_k \geq
2,
$$
which implies inductively that
\be \label {sjk}
    v_{j_k} \geq \sum_{n=0}^{j_k-2} v_n \quad \text{for}\quad j_k \geq 2.
\ee
Letting $k=\kappa$ (so that $j_\kappa= N$ and $v_N=p$) we obtain:
$$
   t: = p -\sum_{n=0}^{N-2} v_n \geq 0.
$$
Since $f^{v_n}(\partial I_{n+1}) \subset \partial  I_n$, 
it follows that $f^p (\partial I_{N-1}) \subset f^t(\partial I_0) \subset \partial I_0$.
\end{pf}

\comm{
\subsubsection{Critical pullback control}

For $0 \leq n \leq N$ and $l \geq m \geq 0$ are such that
$f^{l-m}(0) \in I_{n,m}$, let
$G_{n,l,m}=A \circ f^{l-m} \circ B$ where $A:I_{n,m} \to [-1,1]$ and $B:[-1,1]
\to I_{n,l}(0)$ are orientation preserving affine maps.  We also denote
$G_{n,l}=G_{n,l,0}$.

\begin{thm}

For every $0<\underline \lambda<1$, $\kappa_0 \geq 0$, if $\lambda \geq
\underline \lambda$ and $\kappa \geq \kappa_0$ then
$G_{j_{\kappa-\kappa_0},p}$ belong to a
compact subset $\KK=\KK(\underline \lambda,\kappa_0,\epsilon_0)
\subset C^\omega([-1,1],[-1,1])$.

\end{thm}

\begin{pf}

Let $n=j_{\kappa-\kappa_0}$.

Let $0=l_0<...<l_r=p$ denote the moments such that $0 \in
I_{n,l-l_i}(f^{l_i}(0))$.  Then $G_{n,l}=G_r \circ ... \circ G_1$, where
$G_i=G_{n,l-l_{i-1},l-l_i}$.

Let us say that $G_i,...,$ central cascade is a maximal sequence $G_$
}

\subsection{Transition maps}

The geometric considerations made below are all contained
in \cite {attractors}, though we do not need the finest part of that
argument, dealing with
growth of geometry for Fibonacci-like cascades (with or without saddle-node
subcascades), which is  not valid in higher degree anyway.

We say that $n \geq 0$ is {\it admissible} if $f^n(0) \in I_0$.  For admissible
$n$, let $T_n$ be the closure of the
connected component of $f^{-n}(\inter I_0)$ containing $0$.
In particular $T_0=I_0$.  More generally, letting $w_n=\sum_{k=0}^{n-1} v_k$,
we obtain $T_{w_k}=I_k$.
If $n$ is admissible we let $A_n: T_n \to \I$ be the
orientation preserving affine homeomorphism, where $\I=[-1,1]$.

We say that $T_n$ is a
{\it pullback} of $T_m$ if $n>m$ and
$f^{n-m}(0) \in T_m$.  In this
case, $f^{n-m}$ restricts to a 
map $(T_n, \di T_n) \to (T_m, \di T_m)$, and we let
$G_{n,m}=A_m \circ f^{n-m} \circ A^{-1}_n$, which we call a {\it transition
map}.

We say that $T_n$ is a {\it kid}
of $T_m$ if $T_n$ is a pullback of $T_m$ but is not a pullback of any
$T_k$ with $m<k<n$.  Notice that in this case,
$f^{n-m-1}|f(T_n)$ extends to an analytic diffeomorphism onto $T_m$.
If $T_n$ is a kid of $T_m$, the transition map $G_{n,m}$ is called
{\it short}, otherwise it is called {\it long}.

A short transition map $G_{n,m}$ is called
{\it $\delta$-good} if $f^{n-m-1}|f(T_n)$ extends to an analytic
diffeomorphism onto a $\delta |T_m|$-neighborhood of $T_m$.  The usual
Koebe space argument (see \cite{MvS}) yields:

\begin{lemma} \label {transition}

For every $\delta>0$, any 
$\delta$-good transition map belongs to a compact set
$\KK=\KK(\delta,\epsilon_0) \subset
C^\omega(\I,\I)$, only depending on $\epsilon_0$ and $\delta$.
 
\end{lemma}

Here $C^\omega$ stands for the space of analytic maps, with the usual
inductive limit topology.



There is a unique {\it canonical decomposition} of
a long transition map into short transition maps:
letting $m=n_1<...<n_l=n$
be the sequence of moments such that $f^{n-n_j}(0) \in T_{n_j}$,
then $T_{n_{j+1}}$
is a kid of $T_{n_j}$ for $0 \leq j \leq l-1$ and
$G_{n,m}=G_{n_2,n_1} \circ \cdots \circ G_{n_l,n_{l-1}}$.

A {\it central cascade} is a sequence $T_{n_1},...,T_{n_l}$ such that
$T_{n_{j+1}}$ is the first kid of $T_{n_j}$ and $f^{n_{j+1}-n_j}(0) \in
T_{n_{j+1}}$ for $1 \leq j \leq l-1$.  Notice that in this case
$n_{j+1}-n_j$ is independent of  $j \in [1, l-1]$.
If $n_2-n_1<p$ then we distinguish the {\it saddle-node}
and {\it Ulam-Neumann} types of central cascades according to whether
$0 \notin f^{n_2-n_1}(T_{n_2})$ or $0 \in f^{n_2-n_1}(T_{n_2})$.

A long transition map is called saddle-node/Ulam-Neumann if its canonical
decomposition $G_{n_2,n_1} \circ \cdots \circ G_{n_l,n_{l-1}}$ is such that
$T_{n_1},...,T_{n_l}$ is a saddle-node/Ulam-Neumann cascade.

We say that a long transition map $G_{n_l,n_1}$
is $\delta$-good if all the components  $G_{n_{j+1},n_j}$ of its canonical decomposition
are $\delta$-good. 
Notice that if $G_{n_l,n_1}$ is central then this is equivalent to
$\delta$-goodness of the top level $G_{n_2,n_1}$.

\def\bot{\mathrm{{bot}}}

Besides $\delta$-goodness, an important role is also played by two
parameters associated to a long transition map of saddle-node type: the scaling
factors of the top and bottom levels,
\mbox{ $\lambda_\top=|T_{n_2}|/|T_{n_1}|$ and $\lambda_\bot=|T_{n_l}|/|T_{n_{l-1}}|$.}

\begin{lemma} \label {saddle}

For every $0<\underline \lambda<\overline \lambda<1$,
$\delta>0$, any $\delta$-good
long transition maps of saddle-node type with parameters
$\lambda_\bot,\lambda_\top \in [\underline \lambda,\overline \lambda]$
belongs to a compact set $\KK(\underline \lambda,\overline
\lambda,\delta,\epsilon_0) \subset C^\omega(\I,\I)$.

\end{lemma}

\begin{pf}

For such a transition map $G_{n,m}$, let us consider a
maximal saddle-node cascade $T_{n_1},...,T_{n_L}$ such that $m=n_1$ and
$n=n_l$ for some $l \leq L$.

By Lemma \ref {transition}, we only risk losing compactness when $l$,
and hence $L$, is large, which
is related to the presence of a {\it nearly parabolic fixed point} in
$T_{n_2}$ for $F=f^{n_2-n_1}$ (since there is in fact no fixed point, the
terminology means that a parabolic fixed point appears after a small
perturbation of $F$).  In this case we have the basic geometric estimate,
due to Yoccoz:
\be \label {quadsca}
    \frac {|T_{n_i}|} {|T_{n_{i+1}}|}-1 \sim \max \{i,L-i\}^{-2}, \quad 1 \leq i
\leq L-1
\ee
(the implied constants depending on the bounds $\underline \lambda,\overline
\lambda$ on scaling factors).
See \cite {FM}, Section 4.1, for a discussion of almost parabolic dynamics
and the statement of Yoccoz's Lemma.
\comm{
More precisely, let $\Delta_{L-j}$ be the connected
component of $T_{n_j} \setminus T_{n_{j+1}}$, $1 \leq j \leq L-1$,
which is contained in $F(T_{n_2})$.  The $\Delta_j$ are consecutive
intervals such that $F(\Delta_j)=\Delta_{j+1}$, so $F:\cup_{1 \leq j \leq
L-2} \Delta_j \to \cup_{2 \leq j \leq L-1} \Delta_j$ is an {\it almost parabolic
map} in the terminology of \cite {FM}, Section 4.1.  The upper bound
$\lambda_\top,\lambda_\bot \leq \overline \lambda$
allows us to apply Yoccoz's Lemma (\cite
{FM}, page 354) to get $|\Delta_j| \max \{j,L-1-j\}
\sim \sum_{j=1}^{L-2} |\Delta_j|$, $1 \leq j \leq L-2$, which implies
\be \label {quadsca}
    \frac {|T_{n_i}|} {|T_{n_{i+1}}|}-1 \sim \max \{i,L-i\}^{-2}, \quad 1 \leq i
\leq L-1.
\ee
}
 
In particular, either $l$ is bounded (and we are fine) or $L-l$ is bounded. 
Assuming that $l \geq 4$, $F^3(T_{n_l})$ is
contained in a connected component $J$ of $T_{n_{l-3}} \setminus
T_{n_{L-1}}$.  Since Lemma \ref {transition} provides bounds on $F^3|T_{n_l}$,
we just have to show that $F^{l-4}|J$ is under control.
But the map $F^{l-4}$ maps $J$ onto a connected
component of $T_{n_1} \setminus T_{n_{L-l+3}}$, and extends
analytically to a diffeomorphism onto a connected component $J'$
of $T^\delta_{n_1} \setminus
T_{n_{L-l+4}}$, where $T^\delta_{n_1}$ is a $\delta |T_{n_1}|$-neighborhood of
$T_{n_1}$.  By (\ref {quadsca}), $J'$ is a $\delta' |F^{l-4}(J)|$-neighborhood of
$F^{l-4}(J)$ for some $\delta'>0$, so $F^{l-4}|J$ is under Koebe control.
\end{pf}

\subsection{Small scaling factors}

In \cite {attractors}, several combinatorial properties are shown to yield
small scaling factors.  We will need somewhat simpler estimates, which we
will obtain from the following:

\begin{lemma} \label {easy}

For every $\epsilon>0$, there exists
$\delta=\delta(\epsilon,\epsilon_0)>0$ with the following property.  Assume
that the postcritical set intersects a connected component $D$ of the first
landing map to $I_{j_{k+1}}$ such that $D \subset I_{j_k} \setminus
I_{j_k+1}$ and $|D|/|I_{j_k}|<\delta$.  Then
$|I_{j_\kappa+1}|/|I_{j_k}|<\epsilon$.

\end{lemma}

\begin{pf}

We may assume that $\lambda_{j_k}$ is not small and that $\kappa-k$ is
bounded.  Let $r>0$ be minimal such that $f^r(0) \in D$.

Assume first $r<v_{j_{k+1}}$, i.e.,  $\orb 0$ lands in $D$ before landing in $I_{j_{k+1}}$.
Let us consider the first landing 
$f^{v_{j_k}}(0)\in  I_{j_{k+1}-1}\sm I_{j_{k+1}}$ of $\orb 0$ in $I_{j_k}$.
Let 
\begin{equation}\label{D'}
     D'\subset I_{j_{k+1}-1}\sm I_{j_{k+1}}
\end{equation}
be the pullback of $D$ under $f^{r-v_{j_k}}$ containing $f^{v_{j_k}}(0)$.
(In other words, $D'$ is the component of the domain of
the first landing map to $I_{j_{k+1}}$ containing $f^{v_{j_k}}(0)$.)  

The map  $f^{r-v_{j_k}}: D'\ra D$ is a composition of the
transit map $$f^{(j_{k+1}-j_k-1) v_{j_k}}:I_{j_{k+1}-1}\sm I_{j_{k+1}} \ra I_{j_k}\sm I_{j_k+1}$$
and (possibly) several first return maps to $ I_{j_k}\sm I_{j_k+1}$.
It has a Koebe extension $\tl D' \ra  L$, where $L$ is the connected
component of $I_{j_k} \setminus \{z_0\}$ containing $D$ and
$z_0=f^{(j_{k+1}-j_k-1) v_{j_k}}(0) \in I_{j_k+1}$.  This
implies an upper bound on the derivative of the inverse map $D\ra D'$,
unless the distance from $z_0$ to $D$ is small compared to
$|I_{j_k}|$.

In the former case, we readily conclude that
\begin{equation}\label{D' and D}
|D'| \leq C\, |D|.
\end{equation}
On the other hand, we may assume that $|I_{j_{k+1}-1}|\asymp |I_{j_k}|$
(otherwise the result is obvious), so
$|D'|$ is small compared with $|I_{j_{k+1}-1}|$.
If it is, in fact,  small compared to the distance $\de$ to $\di I_{j_{k+1}-1}$
then 
$$
         \lambda_{j_{k+1}}\asymp  (|D'|/ \de)^{1/d}
$$ 
is small as well.  Otherwise $\de$ is small compared to $| I_{j_{k+1}-1}|$
Then the return to $I_{j_{k+1}-1}$ is ``very low'', i.e.,
$\si:=|f^{v_{j_k}}(I_{j_{k+1}})|/|I_{j_{k+1}-1}|$ is small.  
But then $\lambda_{j_{k+1}}\asymp \si^{1/d}$ is small again.

In the latter case, $z_0$ must be close to $\partial I_{j_k+1}$, so
by considering the transit map
$f^{(j_{k+1}-j_k-2)v_{j_k}}:I_{j_{k+1}-2}
\setminus I_{j_{k+1}} \to I_{j_k} \setminus I_{j_k+2}$, we conclude that
the distance $\delta$ from $f^{v_{j_k}}(0)$ to $\partial I_{j_{k+1}-1}$ is
small compared to $|I_{j_k}|$, so we can apply the above analysis of a
very low return to obtain that $\lambda_{j_{k+1}}$ is small.

\medskip
Assume now $r>v_{j_{k+1}}$, i.e.,    $\orb 0$ lands in  $I_{j_{k+1}}$  before landing in  $D$.
Let $s\in (0,r)$ be its last landing moment in $ I_{j_{k+1}}$
before landing in $D$.  Then 
$$
 f^s (0)  \in I_{j_{k+1}}\sm I_{j_{k+1}+1},\quad 
  \mathrm{for \ otherwise} \quad 
    f^{r-s}(0)\in f^{r-s}( I_{j_{k+1}+1}) \subset D,
$$
contradicting the minimality of $r$. 

Let $\De$ be the pullback of $D$ under $f^{r-s}$ containing $f^s(0)$          
(which is the component of the domain of the first return map to $I_{j_{k+1}}$
containing  $f^s(0)$). The transit map $f^{r-s}: \De \ra D$ is the composition of 
the fist return map $f^{v_{j_k}}$ and the transit map from an interval 
$D'$ to $D$, where $D'$ is defined as (\ref{D'}) except that it is centered at $f^{s+v_{j_k
}}(0)$ rather than at $f^{v_{j_k}}(0)$.
By the previous analysis, either the return to $I_{j_{k+1}-1}$ is very low, and
we are done since this implies that $\lambda_{j_{k+1}}$ is small,
or estimate (\ref{D' and D}) holds.
Since $ f^{v_{j_k}}$ is the power map $x^d$, up to bounded distortion,
and the intervals $I_{j_{k+1}}$,  $I_{ j_{k+1} + 1}$ are comparable,  
we conclude that
$|\De|$ is small compared with  $I_{j_{k+1}}$. 
All the more, the component $\De'$ of the first landing map to $I_{j_{k+2}}$ containing $f^s(0)$
 is small compared with  $I_{j_{k+1}}$. We can now start the
procedure over with $k$ replaced by $k+1$ and $D$ replaced by $\Delta'$.


Since $\kappa-k$ is assumed to be bounded, this process must eventually produce a
small scaling factor.
\end{pf}

One important situation in our analysis corresponds to the critical
orbit hitting deep inside a long central cascade.

More precisely, we say that $T_m$ is $k$-{\it deep} ($k \geq 2$)
inside a central
cascade if there is a central cascade
$T_{n_1}$,..., $T_{n_l=m}$ ,$T_{n_{l+1}}=T_m'$,...,$T_{n_L}$ such that
$k \leq l \leq L-k$.  We say that
the critical orbit {\it hits} $T_m$ if there is $r \geq 0$ such that
$f^r(0) \not\in T'_m$,
but $f^{r+n_2-n_1}(0) \in T_m \setminus T'_m$.

\begin{lemma} \label {sn}

For every $\epsilon>0$ there exists $k=k(\epsilon,\epsilon_0)>0$
with the following property.  Assume that the critical orbit hits some
$T_m$ which is $k$-deep inside a central cascade.
Then $|I_{j_\kappa+1}|/|T_m|<\epsilon$.

\end{lemma}

\begin{pf}

It is no loss of generality to assume that $T_{n_1}=I_{j_i}$ for some $i$
(since between any interval $T_n$ and its kid $T_n'$, there must be an
interval of the principal nest).

We may assume that $\lambda_{j_i}$ is not small.
Let $D'$ be the component of the domain of the first landing map to
$I_{j_{i+1}}$ containing $f^{r+n_2-n_1}(0)$.  Then $|D'|/|I_{j_i}|$ is small,
see (\ref {quadsca}).  Let $s<r+n_2-n_1$ be maximal with
$f^s(0) \in I_{j_i} \setminus I_{j_i+1}$.
Let $D$ be the component of the domain of the first landing map to
$I_{j_{i+1}}$ containing $f^s(0)$.
Pulling back the Koebe space, we get $|D|/|I_{j_i}|$ small.  The result
follows from Lemma \ref {easy}.
\end{pf}

A simpler situation involves long Ulam-Neumann cascades:

\begin{lemma} \label {misiu}

For every $\epsilon>0$ there exists $k=k(\epsilon,\epsilon_0)>0$ with the
following property.  Let $T_{n_1},...,T_{n_k}$ be a central cascade  of
Ulam-Neumann type.  Then there exists $I_t \subset T_{n_k}$ such that
$\lambda_t<\epsilon$.

\end{lemma}

\begin{pf}

As in Lemma \ref {sn}, we may assume that $T_{n_1}=I_{j_i}$.  Due to the
long Ulam-Neumann cascade, $\lambda_{j_{i+1}-1}$ is close to $1$, see \cite[Lemma 8.3]{puzzle}.
In particular, the domain $D$
of the first landing map to $I_{j_{i+1}}$ containing $f^{v_{j_i}}(0)$
has lots of Koebe space in $I_{j_i}$.  Pulling back by $f^{v_{j_i}}$, we
conclude that $|I_{j_{i+1}+1}|/|I_{j_i+1}|$ is small, which implies that
$\lambda_{j_{i+1}}$ is small.
\end{pf}

We will also need the following easy criterion.  Let us say that $T_n$ is
$\delta$-safe if the postcritical set does not intersect a $\delta
|T_n|$-neighborhood of $\partial T_n$.  Notice that if $T_n$ is
$\delta$-safe and $T_m$ is a kid of $T_n$ then $G_{m,n}$ is $\delta$-good.

\begin{lemma} \label {safe}

For every $\epsilon>0$, there exists $\delta>0$ such that if $0 \leq k \leq
\kappa$ is such that $I_{j_k}$ is not $\delta$-safe, then
$|I_{j_\kappa+1}|/|I_{j_k}|<\epsilon$.

\end{lemma}

\begin{pf}

Consider first the case of ``postcritical set inside'', i.e., for some
$r>0$, $f^r(0) \in I_{j_k}$ is near $\partial
I_{j_k}$.  Let $D$ be the component of the first landing map to
$I_{j_{k+1}}$ containing $f^r(0)$.  Since it has Koebe space inside
$I_{j_k}$
(the landing map $D \to I_{j_{k+1}}$ extends to a diffeomorphism $I_{j_k}
\supset D' \to I_{j_k}$),
$|D|/|I_{j_k}|$ is small, so the result follows from Lemma \ref {easy}.

Consider now the case of ``postcritical set outside'', i.e.,
for some $r>0$, $f^r(0) \notin I_{j_k}$ is near $\partial
I_{j_k}$.  We may assume that $k>0$ (otherwise applying $f$ once produces
``postcritical set inside'' reducing to the previous case).  We may assume
further that $\lambda_{j_k}$ is not small (otherwise the result is obvious),
and that if $j_k-j_{k-1}$ is large then this is due to a saddle-node cascade
(otherwise we can apply Lemma \ref {misiu}).  Letting $r'=r+(j_k-j_{k-1}-1)
v_{j_{k-1}}$, it follows that $f^{r'}(0)$ is just outside $\partial
I_{j_{k-1}+1}$.  Let $D$
be the connected component of the first landing map to $I_{j_k}$
containing $f^{r'}(0)$.  Since it has Koebe space inside
$I_{j_{k-1}} \setminus I_{j_{k-1}+1}$,
this implies that $|D|/|I_{j_{k-1}}|$ is small, and the result follows from
Lemma \ref {easy}.
\end{pf}

\subsection{Main precompactness}




\begin{lemma} \label {conca}

Let us consider a composition of transition maps $G_{m_2,m_1} \circ
\cdots \circ G_{m_r,m_{r-1}}$, where each $G_{m_{j+1},m_j}$ is either short
or saddle-node.  Assume that $T_{m_1}$ is $\delta$-safe,
$\lambda(T_{m_1})<1-\delta$ and $|T_{m_r}|/|T_{m_1}|>\delta$.
Assume also that whenever $G_{m_{j+1},m_j}$ is saddle-node of length at
least $l$ then $T_{m_{j+1}}$ is not $l$-deep inside a central cascade.
Then there exists
$\delta'=\delta(\epsilon_0,\delta,l,r)>0$ and a compact subset
$\KK=\KK(\epsilon_0,\delta,l,r) \subset C^\omega(\I,\I)$ such that
for $2 \leq j \leq r$ we have
\begin{enumerate}
\item $T_{m_j}$ is $\delta'$-safe,
\item $\lambda(T_{m_j})<1-\delta'$,
\item $G_{m_j,m_{j-1}} \in \KK$.
\end{enumerate}

\end{lemma}

\begin{pf}

Induction reduces considerations to the case $r=2$.

Notice that if we show that $G_{m_2,m_1}$ is in a compact class, it will
follow that $T_{m_2}$ is $\delta'$-safe (any postcritical set near $\partial
T_{m_2}$ would be taken by $f^{m_2-m_1}$ to postcritical set near $\partial
T_{m_1}$.  Moreover, it will also follow that $\lambda(T_{m_2})<1-\delta'$:
the map $f^{m_2-m_1}$ takes $T_{m_2}'$ into a connected component of the
first landing map to $T'_{m_1}$ and any such component must have Koebe space
inside $T_{m_1}$ since $\lambda(T_{m_1})<1-\delta$.

Thus we just have to show that $G_{m_2,m_1}$ is in a compact class.  Notice
that $G_{m_2,m_1}$ is $\delta$-good.  If
$G_{m_2,m_1}$ is short, the conclusion from Lemma \ref {transition}.  If $G_{m_2,m_1}$
is saddle-node, this will follow from Lemma \ref {saddle} once we show that
$\lambda_\top$ and $\lambda_\bot$ are bounded away from $0$ and $1$. 
Clearly both are at least $\delta$ and moreover
$\lambda_\top=\lambda(T_{m_1})<1-\delta$.

Let $T_{n_1=m_1}$,...,$T_{n_s=m_2}$,...,$T_{n_L}$ be the maximal central
cascade starting at $T_{n_1}$.
As in Lemma \ref {saddle}, see (\ref {quadsca}), we see that if $\lambda_\bot$ is close
to $1$ then $s$ and $L-s$ are large.  But by hypothesis $\min \{s,L-s\} \leq
l$, giving the result.
\end{pf}


The following two similar estimates will be proved simultaneously:

\begin{lemma} \label {mapre}

For $\epsilon>0$ 
there exists a compact subset
$\KK=\KK(\epsilon,\epsilon_0) \subset C^\omega(\I,\I)$
with the following property.
Assume that $|I_{j_\kappa+1}|/|I_{j_k}|>\epsilon$.
Then $G_{w_{j_k}+p,w_{j_k}} \in \KK$.

\end{lemma}

\begin{lemma} \label {mapre1}

For $\epsilon>0$ and $b_0 \in \N$,
there exists a compact subset
$\KK=\KK(\epsilon,\epsilon_0,b_0) \subset C^\omega(\I,\I)$
with the following property.
Assume that $|I_{j_\kappa+1}|/|I_0|>\epsilon$.  If $0 \leq b <\min \{p,b_0\}$ is such that
$p-b$ is admissible then
$G_{p-b,0} \in \KK$.

\end{lemma}

\subsubsection{Proof of Lemmas \ref {mapre} and \ref {mapre1}}

The proofs of both lemmas follow a basically parallel path.
In both cases we need to estimate a map of type
$G_{w_{j_k}+p-b,w_{j_k}}$, where in the setting of Lemma \ref
{mapre} we set $b=0$, 
while in the setting of Lemma \ref {mapre1} we set $k=0$.  
Write it as a composition
of a minimal number of
either short transition maps or long transition maps of Ulam-Neumann or
saddle-node type.
Clearly the number of elements of this decomposition is bounded in terms
of $\kappa-k$, which in turn is bounded in terms of $\epsilon$.

Let us now consider a finer decomposition
$G_{m_2,m_1} \circ \cdots \circ G_{m_r,m_{r-1}}$,
where we split the Ulam-Neumann pieces into short
transition maps (so that each $G_{m_{j+1},m_j}$ is either short or
saddle-node).  The Ulam-Neumann cascades have
bounded length (by Lemma \ref {misiu}), so $r$ is also bounded.

In order to conclude, it is enough to show that the conditions of
Lemma \ref {conca} are satisfied.  Since
$T_{m_1}=I_{j_k}$,
$\lambda(T_{m_1})$ is indeed bounded away from $1$, and
it is $\delta$-safe by Lemma \ref {safe}.  So we just have to check that for $1
\leq j \leq r-1$, if $G_{m_{j+1},m_j}$ is saddle-node with big length then
$T_{m_{j+1}}$ is not too deep inside a central cascade.

The following combinatorial estimate
will be key to the analysis.

\begin{lemma} \label {G}

Suppose that $G_{m_{j+1},m_j}$ is saddle-node, and let
$T_{n_1=m_j}$,...,$T_{n_l=m_{j+1}}$ be the associated cascade.  If $l \geq
4+b$ then $f^{m_r-m_j}(0) \notin T_{n_{4+b}}$.

\end{lemma}

\begin{pf}

Let $j_k+1 \leq s \leq j_\kappa$ be minimal such that
$T_{m_j} \supset I_s$.  It follows that 
\be\label{sum1}
m_j-m_1 \leq \sum_{n=j_k}^{s-1} v_n.
\ee

Since $I_s \subset T_{n_1} \subset \inter I_{s-1}$, it follows
that $I_{s+1} \subset T_{n_2} \subset \inter
I_s$ and $T_{n_3} \subset \inter I_{s+1}$.
Then $n_3-n_2 \geq v_s \geq n_2-n_1$, and since $n_{t+1}-n_t=n_2-n_1$ for $1
\leq t \leq l-1$, we see that $n_2-n_1=v_s$.

Assume first that $s=j_t$ for some $k<t \leq \kappa$.  Then 
(\ref{sjk}) and (\ref{sum1}) imply
$m_j-m_1 \leq 2 v_s$, so that $m_j-m_1+b\leq (2+b)v_s$.  
Thus, 
$$
p = m_r-m_1+b= (m_r-m_j)+(m_j-m_1)+b = (m_r-m_j) + q, \quad \text{where}\
    q\leq (2+b)v_s.
$$
If  $x:= f^{m_r-m_j}(0) \in T_{n_{3+b}}$ then $f^{n v_s} (x)\in f^{v_s}
T_{n_2}$ for  $n\leq 2+b$.
Hence
$f^p(0)=f^q(x)$ is either
in $f^{v_s}(T_{n_2})$ or it is outside $T_{n_1}$.  In any case, it can
not belong to the renormalization interval $\cap_{n \geq 1} I_n$ --
contradiction.

Assume now that $j_t+1=s<j_{t+1}$ for some $k \leq t<\kappa$.  Arguing as
before, we see that $m_j-m_1 \leq 3 v_s$, and if $f^{m_r-m_j}(0) \in
T_{n_{4+b}}$ we arrive at a similar contradiction.

Assume now that $j_t+2 \leq s<j_{t+1}$ for some $k \leq t<\kappa$.      
Then the map $f^{n_2-n_1}$ has a unimodal extension to the interval $I_{s-1}\supset T_{n_1}$.
Hence $T_{n_1}$ is a kid of the interval $T_{n_0}\subset I_{s-2}$ of depth $n_0=2n_1-n_2$.  
But then  $G_{m_{j+1},m_j}$ is not a maximal saddle-node transition map in the
decomposition of $G_{m_r,m_1}$, contradicting the definition of $m_j$.
\end{pf}

Let $T_{n_1=m_j}$,...,$T_{n_l=m_{j+1}}$,...,$T_{n_L}$ be the maximal
continuation of the saddle-node cascade associated to $G_{m_{j+1},m_j}$, and
assume that $l$ and $L-l$ are large.
By Lemma~\ref{G}, $f^{m_r-m_j}(0)
\notin T_{n_{4+b}}$, which implies that $f^{m_r-m_{j+1}}(0) \in T_{n_s} \setminus
T_{n_{s+1}}$ for some $l \leq s \leq l+2+b$.  
Indeed,  $f^{m_r-m_{j+1}}(0) \in T_{m_{j+1}}= T_{n_l}$,
but   $f^{m_r-m_{j+1}}(0) \not\in T_{n_{l+3+b}}$, 
for otherwise 
$$
    f^{m_r-m_j}(0) = f^{m_{j+1}-m_j}(f^{m_r-m_{j+1}}(0))\in f^{n_l-n_1} ( T_{n_{l+3+b}} ) \subset T_{n_{b+4}}.  
$$

On the other hand, since
$G_{m_{j+1},m_j}$ is a maximal saddle-node cascade in the decomposition of
$G_{m_r,m_1}$, we must have $f^{m_r-m_{j+1}+n_1-n_2}(0) \notin
T_{n_{s+1}}$.  
We can then apply Lemma \ref {sn} to conclude that
$|I_{j_\kappa+1}|/|I_{j_k}|$ is small, contradiction.
This establishes that either $l$ or
$L-l$ must be small, as desired.
\qed



\subsection{Proof of Lemma \ref {Dichotomy}}





We may assume that the sequence $f_n=\tilde f_n$ converges to some $f_\infty$.
Let $p_n$ be the period of $f_n$.
Let $\Lambda_n$ be the affine map such that $\Lambda_n \circ (f_n)^{p_n} \circ
\Lambda_n^{-1}$ is normalized.

Assume first that the prerenormalization intervals of $f_n$
do not have length
bounded from below: following
the terminology of \cite {puzzle}, we will say that
the combinatorics of the $f_n$ is not essentially
bounded.  Then either
$\inf \lambda_{N(f_n)}(f_n)=0$ or $\sup \kappa(f_n)=\infty$ by \cite {attractors}.

If $\inf \lambda_{N(f_n)}(f_n)=0$ then $f_\infty$ is a unicritical polynomial \cite
{attractors}.

Consider now the case $\inf \lambda_{N(f_n)}(f_n)>0$ and $\sup
\kappa(f_n)=\infty$.  We may assume
that $\lim \kappa(f_n)=\infty$. 
Passing through
a subsequence we may assume that for each $k \geq 0$,
$\Lambda_n(I_{j_{\kappa(f_n)-k}}(f_n))$
converges to a closed interval $D_k$.
Clearly each $D_k$ is a bounded interval (scaling factors minorated)
and $\cup D_k=\R$ (scaling factors bounded away from $1$).
We may also assume that $\Lambda_n(T_{w_{j_{\kappa(f_n)-k}}+p_n}(f_n))$ converges to a closed
interval $D'_k$.  Then $D'_k \subset D_1$ by Lemma \ref {I'_0}.
By Lemma \ref {mapre}, for every $k \geq 0$, $f_\infty$ has an analytic extension
$D'_k \to D_k$ which is proper.  It follows that $f_\infty$ has a maximal
analytical extension to $\cup_{k \geq 0} D'_k \subset D_1$.

Assume now that $f_n$ has essentially bounded, but unbounded,    
combinatorics.  We may assume that $p_n \to \infty$.  Let
$0=b^0_n<b^1_n<...$ be the sequence of admissible moments, i.e., such that
$f_n^{b^i_n}(0) \in I_0(f_n)$.  
Clearly $b^i_n \leq i b^1_n$ (since the critical point returns to $I_0$ no earlier than any other point $x\in I_0$). 
Moreover, $b^1_n$ is bounded (otherwise the combinatorics is close to the
Chebyshev one and
we would already have $\inf |I_1(f_n)|=0$).  Notice that if $0<b^i_n<p_n$,  
$f_n^{p_n}(\partial T_{p_n-b^i_n})$ is the orientation reversing fixed point of $f_n$.
Hence $f_n^{p_n}$ has at least $2i-1$ critical points in $T_{p_n-b^i_n}(f_n)$
 (counted with multiplicity).

We may assume that the intervals $\Lambda_n(T_{p_n-b^i_n})$
converge to intervals $D_i$ for each $i$.
Clearly $\cup D_i$ is a bounded interval.
By Lemma \ref {mapre1}, $f_\infty$ has an analytic extension to $\cup D_i$, and
restricted to each $D_i$ it has at least $2i-1$ critical points.  
So $f_\infty$ cannot extend beyond $\cup D_i$.
\qed

\comm{
\subsection{Proof of Theorem \ref {limits of towers}}

It is enough to consider the case $n=-1$ (otherwise we just apply induction
and shift the indices).

Notice that $f^{(\infty)}_0$ is affinely conjugated to an appropriate
restriction of $(f^{(\infty)}_n)^{t_n}$, where $t_n=p_n \cdots p_{-1}$ and
$p_i$ is the period of $f^{(\infty)_i}$.  In particular, it has an analytic
extension, with several critical points, to arbitrarily large intervals.
By Theorem \ref {singularity},
the renormalization combinatorics of $f^{(k)}_{-1}$ is
bounded.  The result now immediately follows from:

\begin{lemma} \label {esbound}

Let $f^{(k)} \in \CC^\R(\epsilon)$ be a sequence of renormalizable maps
with bounded combinatorics.  If $\lim R f^{(k)}$ is the
renormalization of some $f \in \CC^\R$, then the combinatorics of
$f^{(k)}$ is eventually the same as that of $f$ and $f^{(k)} \to f$.

\end{lemma}

\begin{pf}

The set of renormalizable maps in $\CC^\R$ with a given combinatorics is
closed, and the restriction of the renormalization operator to this closed
set is continuous.  Thus all limits of $f^{(k)}$ are renormalizable, and
their renormalization must be $R f$.
Since the renormalization operator is injective for real p-l maps
\cite {MvS}, $f^{(k)}$ may only accumulate on $f$.
\end{pf}
}

\end{document}